\documentclass[fleqn]{elsarticle}
\usepackage{amssymb}
\usepackage[sub,ovp]{psfragx}
\usepackage{overpic,color}
\usepackage{psfrag}
\usepackage{natbib}
\usepackage{amsbsy}
\usepackage{upgreek}
\usepackage{wasysym}	
\DeclareSymbolFont{UPM}{U}{eur}{m}{n}
\DeclareMathSymbol{\uppartial}{0}{UPM}{"40}

\def\bs{\boldsymbol{\Sigma}}
\def\bS{\boldsymbol{\Sigma}}
\def\bO{\boldsymbol{\Omega}}
\begin{document}
\begin{frontmatter}
\title{ An Exact Consistent Tangent Stiffness Matrix for a Second Gradient Model for Porous Plastic Solids: Derivation and Assessment. }
\author{Koffi Enakoutsa$^{1,2}$ }
\address{$^1$Department of Mathematics, California State University, Northridge, 18111 Nordhoff St, Northridge, CA 91330}
\address{$^2$ Department of Mathematics, University of California Los Angeles, 520 Portola Plaza, Los Angeles, CA 90095, email: koffi@math.ucla.edu}
\begin{abstract}
It is well known that the use of a consistent tangent stiffness matrix is critical to obtain quadratic convergence of the global Newton iterations in the finite element simulations of problems involving elasto-plastic deformation of metals, especially for large scale metallic structure problems.  
In this paper we derive an exact consistent stiffness matrix for a porous material model, the GLPD model developed by Gologanu, Leblond, Perrin, and Devaux for ductile fracture for porous metals based on generalized continuum mechanics assumptions.
Full expressions for the derivatives of the Cauchy stress tensor and the generalized moments stress tensor the model involved are provided. 
The effectiveness and robustness of the proposed tangent stifness moduli are assessed by applyting the formulation in the finite element simulations of ductile fracture problems. 
Compraisons between the performance our stiffness matrix and the standard ones are also provided.         
\end{abstract}
\end{frontmatter}
\newpage
\section{Introduction}
\label{sec:Intro}
Constitutive models involving softening all predict unlimited localization of strain and damage.
This feature generates such undesired phenomena as absence of energy dissipation during crack propagation and mesh size sensitivity in
finite element computations.
Gurson \cite{G77}'s famous model for porous ductile materials, which was derived from approximate limit-analysis of some elementary voided cell in a plastic solid, is no exception.
In this model, unlimited localization arises from the softening because of the gradual increase of the porosity.
$\\$   

Several proposals have been made to solve this problem. One of these, due to Leblond {\it et al.} \textcolor{blue} {\cite{Leb94}} but based on a previous suggestion made by Pijaudier {\it et al.} \textcolor{blue} {\cite{PijB87}} in damage of concrete, comprises adopting a nonlocal evolution equation for the porosity involving some spatial convolution of some ``local porosity rate'' within an otherwise unmodified Gurson model. This simple proposal has attracted the attention of several authors (Tvergaard and Needleman \textcolor{blue} {\cite{TN95}}, Tvergaard and Needleman \textcolor{blue} {\cite{TN97}}, Enakoutsa {\it et al.} \textcolor{blue} {\cite{E07, ELP07}}). It was notably checked by Tvergaard and Needleman \textcolor{blue} { \cite{TN95}} that it allows to eliminate mesh size effects. Also, Enakoutsa {\it et al.} \textcolor{blue} {\cite{E07, ELP07}} showed that with a minor modification, it leads to great numerical reproduction of the results of typical experiments of ductile rupture.
$\\$ 

One shortcoming of Leblond {\it et al.} \textcolor{blue} { \cite{Leb94}}'s proposal, however, is that it is purely heuristic and lacks any serious theoretical justification. This was the motivation for a later, more elaborate and physically based proposal of Gologanu {\it et al.} \textcolor{blue} { \cite{GLPD97}}. These authors derived an improved variant of Gurson's model (the {\it GLPD model}\footnote{GLPD: Gologanu-Leblond-Perrin-Devaux.}) through some refinement of this author's original homogenization procedure based on Mandel \textcolor{blue} { \cite{M64}}'s and Hill \textcolor{blue} { \cite{H67}}'s classical conditions of homogeneous boundary strain rate. In the approach of Gologanu {\it et al.} \textcolor{blue} { \cite{GLPD97}}, the boundary velocity is assumed to be a quadratic, rather than linear, function of the coordinates. The physical idea is to account in this way for the possibility of quick variations of the macroscopic strain rate, such as encountered during strain localization, over short distances of the order of the size of the elementary cell considered. The output of the homogenization procedure is a model of ``micromorphic'' nature involving the second gradient of the macroscopic velocity and generalized macroscopic stresses of ``moment'' type (homogeneous to the product of a stress and a distance), together with some ``microstructural distance'' connected to the mean spacing between neighboring voids.
$\\$

The numerical implementation of the GLPD model into a finite element code is quite involved as this task required to introduce extra degrees of freedom representing strains; these extra degrees of freedom will permit the calculation of the spacial derivative of the strains, but their number will increase from 2 to 6 in 2-dimensional calculations and from 6 to 9 in the 3-dimensional calculations. This will increase the CPU time for the simulations. Another difficulty lies in the necessary operation of “projection” onto the sophisticated yield locus. An implicit algorithm similar in principle to that classically used for the von Mises criterion, although much more complex in detail, is adopted for this purpose. Convergence of the global elastoplastic iterations was difficult. 
To preserve the quadractic convergence rate of the global Newton iterations stifness tangent moduli are needed. The derivation of tangent moduli is a very difficult task, especially for constitutive models with complex forms.  
In this paper we derive an exact consistent stiffness matrix for the GLPD model for porous materals. The derivation is based on the small strain formulation, this choice is consistent with the one adopted in many finite element codes, including Abaqus$^{\circledR}$, LS-Dyna$^{\circledR}$, Adina$^{\circledR}$, and Systus$^{\circledR}$ which have demonstrated their efficiency.
In Section 2 we summarize the consitutive equations of the GLPD model. Section 3 presents some aspects of the numerical implementation of the GLPD model into finite element codes. 
In Section 4 we provide the expressions for the derivatives of the Cauchy stress tensors and the generalized moments stress tensor the model involves with respect to the main field variables. We do not compute all the terms of the tangent stiffness matrix, but only the ones that are critical for the numerical implementation. 
Finally, Section 5 assesses the performance of the tangent stiffness matrix toward its capacity to reach quadratic convergence in the simulations of small scale ductile fracture problems.       
\section{The GLPD Model}
\label{sec: GLPDModel}
The purpose of this section is to provide a complete description of the GLPD model developed by Gologanu, Leblond, Perrin, and Devaux. The original reference \textcolor{blue}{\cite{GLPD97}} for the GLPD model is not easily accessible, a summary of the equations of this model is given here.
A short presentation of the derivation of these equations derived from some homogenization procedure is also provided below; strictly
speaking, this presentation is not indispensable, but it is useful to fully grasp the physical foundations of the GLPD model.

\subsection{Generalities}
\label{subsec:Gen}
In the GLPD model, internal forces are represented through some ordinary second-rank symmetric Cauchy stress tensor ${\bS}$ plus some additional
third-rank ``moment tensor'' ${\bf M}$ symmetric in its first two indices only\footnote{The component $M_{ijk}$ is noted $M_{k|ij}$ in \textcolor{blue}{\cite{GLPD97}'s}
original paper.
The present notation leads to more natural-looking expressions.}.
%
%(The components of ${\bs}$ and ${\bf M}$ are interpreted in paper \cite{GLPD97} as the average values and ``moments'' of the
%components of the microscopic stress tensor over the elementary cell considered, but this interpretation will play no role here).
%
The components of ${\bf M}$ are related through the three conditions
\begin{equation}\label{eqn:RelMoments}
    M_{ijj} = 0.
\end{equation}
(These conditions may be compared to the condition of plane stress in the theory of thin plates or shells).

The virtual power of internal forces is given by the expression
\begin{equation}\label{eqn:VirtPowIntForces}
    {\mathcal P}^{(i)} \equiv - \int_{\Omega} ( {\bS}:{\bf D} + {\bf M}\,\vdots\,{\nabla}{\bf D} ) \, d\Omega
\end{equation}
where $\Omega$ denotes the domain considered, ${\bf D}\equiv \frac{1}{2}\left[ \nabla{\bf V} + (\nabla{\bf V})^T \right]$ (${\bf V}$: material velocity)
the Eulerian strain rate, ${\nabla}{\bf D}$ its gradient, ${\bS}:{\bf D}$ the double inner product $\Sigma_{ij}D_{ij}$ and ${\bf M}\,\vdots\,{\nabla}{\bf
D}$ the triple inner product $M_{ijk}D_{ij,k}$.
$\\$

The virtual power of external forces is given by
\begin{equation}
\label{eqn:VirtPowExtForces}
    {\mathcal P}^{(e)} \equiv \int_{\partial \Omega} {\bf T}.{\bf V} \, dS
\end{equation}
where ${\bf T}$ represents some surface traction\footnote{The general equilibrium equations and boundary conditions corresponding to the expressions
\textcolor{blue}{(\ref{eqn:VirtPowIntForces})} and \textcolor{blue}{(\ref{eqn:VirtPowExtForces})} of the virtual powers of internal and external forces need not be given since they are not
necessary for the numerical implementation.}.
$\\$

The corresponding equilibrium equations read, in the absence of body forces and moments:
\begin{equation}
\label{eqn:Equil}
\Sigma_{ij, j} - M_{ijk, jk} = 0 \quad \mbox{in } \Omega.
\end{equation}
The boundary conditions are complex and will not be given here.
(In fact the numerical implementation of the model will require neither the equilibrium equations nor the boundary conditions, but the
sole expression of the virtual power of internal forces).
$\\$

The hypothesis of additivity of elastic and plastic strain rates reads
\begin{equation}
\label{eqn:AddStrainRate}
    \left\{
        \begin{array}{lll}
          {\bf D} & \equiv & {\bf D}^e + {\bf D}^p \\
          \nabla {\bf D} & \equiv & (\nabla{\bf D})^e + (\nabla{\bf D})^p.\\
        \end{array}
    \right.
\end{equation}
The elastic and plastic parts $(\nabla{\bf D})^e$, $(\nabla{\bf D})^p$ of the gradient of the strain rate here do {\it not} coincide in general with the
gradients $\nabla({\bf D}^e)$, $\nabla({\bf D}^p)$ of the elastic and plastic parts of the strain rate.
\subsection{Hypoelasticity law}
\label{subsec:Hypo}
The elastic parts of the strain rate and its gradient are related to the rates of the stress and moment tensors through the following hypoelasticity law:
\begin{equation}
\label{eqn:Hypoelast}
    \left\{
        \begin{array}{lll}
          \displaystyle \frac{D\Sigma_{ij}}{Dt} & = & \lambda\, \delta_{ij} D_{kk}^e  + 2\mu D_{ij}^e \\
          \displaystyle \frac{D M_{ijk}}{Dt} & = & \displaystyle \frac{b^2}{5} \biggl[ \lambda\, \delta_{ij} (\nabla D)_{hhk}^e  + 2\mu (\nabla D)_{ijk}^e\\
&& -2 \lambda\, \delta_{ij} U_k^e  - 2\mu \left( \delta_{ik} U_j^e  + \delta_{jk} U_i^e \right) \biggr].\\
        \end{array}
    \right.
\end{equation}
In these expressions $\lambda$ and $\mu$ denote the Lam\'e coefficients and $b$ the mean half-spacing between neighboring voids.
(In the homogenization procedure, $b$ is the radius of the spherical elementary cell considered).
Also, $\frac{D \Sigma_{ij}}{Dt}$ and $\frac{ D M_{ijk}}{D t}$ are the Jaumann (objective) time-derivatives of $\Sigma_{ij}$ and ${M}_{ijk}$, given by
\begin{equation}\label{eqn:JaumannDeriv}
    \left\{
        \begin{array}{lll}
          \displaystyle \frac{D \Sigma_{ij}}{D t} & \equiv & \dot{ \Sigma}_{ij} + \Omega_{ki}\Sigma_{kj} + \Omega_{kj}\Sigma_{ik}\\
          \displaystyle \frac{D M_{ijk}}{D t} & \equiv & \dot{M}_{ijk} + \Omega_{hi}M_{hjk} + \Omega_{hj}M_{ihk} + \Omega_{hk}M_{ijh}\\
        \end{array}
    \right.
\end{equation}
where ${\bO} \equiv \frac{1}{2}\left[ {\nabla}{\bf V} - ({\nabla}{\bf V})^T \right]$ is the antisymmetric part of the velocity gradient.
Finally ${\bf U}^e$ is a vector the value of which is fixed by equations \textcolor{blue}{(\ref{eqn:RelMoments})} (written in rate form, $\frac{D{M}_{ijj}}{Dt}=0$):
\begin{equation}\label{eqn:Ue}
    U_i^e = \frac{ \lambda (\nabla D)_{hhi}^e + 2\mu (\nabla D)_{ihh}^e}{ 2\lambda + 8\mu }.
\end{equation}
(This vector may be compared to the through-the-thickness component of the elastic strain rate in the theory of thin plates or shells, the value of which
is fixed by the condition of plane stress).
\subsection{Yield criterion}
\label{subsec:Yield}
The plastic behavior is governed by the following Gurson-like criterion:
\begin{equation}
\label{eqn:Crit}
    \frac{1}{{\Sigma}^2} \left( \Sigma_{eq}^2 + \frac{Q^2}{b^2} \right)
                            + 2 p \, \cosh \left( \frac{3}{2} \frac{\Sigma_m}{{\Sigma}} \right) - 1 - p^2 \leq 0.
\end{equation}
In this expression:
\begin{itemize}
    \item $\Sigma_{eq} \equiv \left( \frac{3}{2}{\bS}':{\bS}' \right)^{1/2}$ (${\bS}'$: deviator of ${\bS}$) is the
          von Mises equivalent stress.
    \item $\Sigma_m \equiv \frac{1}{3}\,{\rm tr}\,{\bS}$ is the mean stress.
    \item ${\Sigma}$ represents a kind of average value of the yield stress in the heterogeneous metallic matrix, the
          evolution equation of which is given below.
    \item $p$ is a parameter connected to the porosity (void volume fraction) $f$ through the relation:
\begin{equation}
\label{eqn:p}
    p \equiv qf^*,f^* \equiv
    \left\{
       \begin{array}{lll}
         f & \mbox{if $f \leq f_c$} \\
         f_c + \delta (f-f_c) & \mbox{if  $f > f_c$} \\
       \end{array}
    \right.
\end{equation}
where $q$ is {\it Tvergaard's parameter}, $f_c$ the {\it critical porosity at the onset of coalescence of voids}, and $\delta$ ($>1$) a factor describing the accelerated degradation of the material during coalescence \textcolor{blue}{\cite{T81, TN84}},
    \item $Q^2$ is a quadratic form of the components of the moment tensor given by
\begin{equation}\label{eqn:Q}
    Q^2 \equiv A_{I} M_{I} + A_{II} M_{II} \quad , \quad
    \left\{
       \begin{array}{lll}
         A_{I} & = & 0.194 \\
         A_{II} & = & 6.108 \\
       \end{array}
    \right.
\end{equation}
         where $M_{I}$ and $M_{II}$ are the quadratic invariants of ${\bf M}$ defined by:
\begin{equation}\label{eqn:M1M2}
  \left\{
    \begin{array}{lll}
      M_{I} & \equiv &  M_{mi} M_{mi} \\
      M_{II} & \equiv & \frac{3}{2} M'_{ijk} M'_{ijk},
    \end{array}
  \right.
\end{equation}
$M_{mi} \equiv \frac{1}{3}M_{hhi}$ and ${\bf M}'$ denoting the mean and deviatoric parts of ${\bf M}$, taken over its first two indices.
   \item Again, $b$ is the mean half-spacing between neighboring voids.
\end{itemize}
\subsection{Flow rule}
\label{subsec:FlowRule}
The plastic parts of the strain rate and its gradient are given by the flow rule associated to the criterion \textcolor{blue}{(\ref{eqn:Crit})} through normality:
\begin{equation}\label{eqn:FlowRule}
\left\{\begin{array}{lll} 
D_{ij}^p & = & \displaystyle H \frac{\partial \Phi}{\partial \Sigma_{ij}}({\bS},{\bf M},\Sigma,f) \\
(\nabla D)_{ijk}^p & = & \displaystyle H \frac{\partial \Phi}{\partial M_{ijk}}({\bS},{\bf M},\Sigma,f)+ \delta_{ik} U_j^p \\
&& ~~~~~~~~~~~~~+ \delta_{jk} U_i^p \end{array}
    \right\}
\end{equation}
were as $$H = \left\{\begin{array}{lll} 
= 0 & \quad {\rm if} \quad & \Phi ({\bS},{\bf M},\Sigma,f) < 0 \\
\geq 0 & \quad {\rm if} \quad & \Phi ({\bS},{\bf M},\Sigma,f) = 0 \end{array} \right\}$$

The term $\delta_{ik} U_j^p + \delta_{jk} U_i^p$ in equation \textcolor{blue}{(\ref{eqn:FlowRule})} represents a rigid-body motion of the elementary cell, which is
left unspecified by the flow rule but fixed in practice by conditions \textcolor{blue}{(\ref{eqn:RelMoments})}.
(The vector ${\bf U}^p$ may be compared to the through-the-thickness component of the plastic strain rate in the theory of thin plates or shells, the
value of which is fixed by the condition of plane stress).
$\\$

The values of the derivatives of the yield function $\Phi({\bS},{\bf M},\Sigma,f)$ in equations \textcolor{blue}{(\ref{eqn:FlowRule})} are readily calculated to be
\begin{equation}
\label{eqn:DerivYieldFunc}
\left\{ \begin{array}{lll} \displaystyle \frac{\partial \Phi}{\partial \Sigma_{ij}}({\bS},{\bf M},\Sigma,f) & = & \displaystyle 3\frac{\Sigma'_{ij}}{{\Sigma}^2} + \frac{p}{\Sigma} \,  \delta_{ij} \sinh \left( \frac{3}{2} \frac{\Sigma_m}{\Sigma} \right) \\
\displaystyle \frac{\partial \Phi}{\partial M_{ijk}}({\bS},{\bf M},\Sigma,f) & = & \displaystyle \frac{1}{{\Sigma}^2 b^2}\left( \frac{2}{3} A_1 \delta_{ij} M_{mk}  + 3 A_2 M'_{ijk}\right) \end{array} \right\}
\end{equation}
\subsection{Evolution of internal parameters}
\label{subsec:EvolIntParam}
The evolution of the porosity is governed by the classical equation resulting from approximate incompressibility of the metallic matrix:
\begin{equation}\label{eqn:EvolPoro}
    \dot{f} = (1-f)\, {\rm tr}\,{\bf D}^p.
\end{equation}
The parameter ${\Sigma}$ is given by
\begin{equation}
\label{eqn:Sigbar}
    {\Sigma} \equiv \Sigma(E)
\end{equation}
where $\Sigma({\epsilon})$ is the function which provides the yield stress of the matrix material in terms of the local equivalent cumulated plastic
strain $\epsilon$, and $E$ represents some average value of this equivalent strain in the heterogeneous matrix.
The evolution of $E$ is governed by the following equation:
\begin{equation}\label{eqn:EvolEpsBar}
    (1-f) {\Sigma} \dot{E} = {\bS}:{\bf D}^p + {\bf M}\,\vdots\, (\nabla{\bf D})^p.
\end{equation}
\section{Numerical implementation}
\label{sec:Impl}
The GLPD model described in Section \ref{sec: GLPDModel} has been incorporated into the Systus$^{\circledR}$ FE code developed by ESI Group, in the 2D case.
The trickiest features of the numerical implementation, which stands as an extension of those proposed by Aravas  \textcolor{blue} { \cite{A87}} and Enakoutsa {\it et al.}  \textcolor{blue} { \cite{ELP07}} for the original
Gurson model, are presented here.
Emphasis is mainly placed on the complex problem of projection of the (supposedly known) elastic stress predictor onto the yield locus defined by the
yield function (\ref{eqn:Crit}). (This problem will be called the {\it projection problem} for shortness in the sequel).
\subsubsection{The GLPD model and the class of generalized standard materials}
\label{subsec:StandGenMater}
The class of {\it generalized standard materials}, as defined by Halphen and Nguyen  \textcolor{blue} {\cite{HN75}}, consists of elastic-plastic materials for which the plastic strain {\it
plus} the internal parameters collectively obey some ``extended normality rule''.
This class is remarkable in that as shown by Nguyen  \textcolor{blue} {\cite{N77}}, for such materials, provided that the flow rule is discretized in time with an {\it implicit}
(backward Euler) scheme, the projection problem is equivalent to minimizing some strictly convex function, which warrants existence and uniqueness of its
solution.
$\\$

It so happens that the GLPD model fits into the framework of generalized standard materials {\it for a fixed porosity}. This property is tied to the special evolution equation (\ref{eqn:EvolEpsBar}) obeyed by the hardening parameter $E$. The proof is provided in Enakoutsa's  \textcolor{blue} { \cite{E07}}'s thesis and is in fact a straightforward extension of that given by Enakoutsa {\it et al.}  \textcolor{blue} {\cite{ELP07}} for the original Gurson model.
$\\$

This property strongly suggests adopting an {\it implicit} algorithm to solve the projection problem, to take advantage of the guaranteed existence and uniqueness of the solution. However, since the porosity $f$ must not be allowed to vary for the GLPD model to be ``generalized standard'', it appears necessary, to benefit from this property, to use an {\it explicit} scheme regarding this specific parameter. Then $f$ will be fixed during the whole calculation of the values of field quantities at time $ t + \Delta t$ from their values at time $t$, and updated (using a discretized version of equation (\ref{eqn:EvolPoro})) only at the end upon the convergence; the projection algorithm will then be exactly the same {\it as if} the porosity were a constant.
$\\$

We shall therefore use the explicit estimate of the porosity at time $t+\Delta t$ given by
\begin{equation}
\label{eqn:ApproxPoroT+DT}
  f ( t + \Delta t) \simeq f(t) + \dot{f} (t) \Delta t,
\end{equation}
and the explicit estimate of the parameter $p(t+\Delta t)$ resulting from there, during the whole ``transition from time $t$ to time $t + \Delta t$'', but
the projection algorithm developed will otherwise be fully implicit with respect to all other parameters, that is the components of the plastic strain
and the plastic strain gradient and the hardening parameter $E$.
From now on, all quantities will conventionally be denoted with a lower index $_0$ if considered at time $t$, and without any special symbol if
considered at time $t + \Delta t$. From now on, all quantities will implicitly be considered at time $t + \Delta t$.
% and denoted without any particular symbol.
%%
\subsubsection{Parametrization of the yield locus}
\label{subsec:ParamYieldLoc}
One key point of the procedure of solution of the projection problem, aimed at reducing the number of unknowns, lies in a suitable partial
parametrization of the yield locus defined by the yield function (\ref{eqn:Crit}).
This parametrization is inspired from the classical one for an ellipse and obtained by looking for the maximum possible value of the quantity
$\Sigma_{eq}^2 + {Q^2}/{b^2}$, namely $(1-p)^2 {\Sigma}^2$, and then writing this quantity in the form $(1-p)^2 \, {\Sigma}^2 \cos ^2 \phi$ for some
angle $\phi$ and solving the equation $\Phi ({\bS},{\bf M},\Sigma,f)=0$ with respect to $\Sigma_m$.
One thus gets
\begin{equation}
\label{eqn:ParamCrit}
    \left\{
        \begin{array}{lll}
            \displaystyle \Sigma_{eq}^2 + \frac{Q^2}{b^2} & \equiv & (1-p)^2 \, {\Sigma}^2 \cos ^2 \phi \\
            \Sigma_m & \equiv & \displaystyle \frac{2}{3} \, {\Sigma} \ \mbox{sgn}(\phi) \ {\rm arg}\,{\rm cosh}
                                \left[ 1 + \frac{ (1-p)^2 \sin ^2 \phi }{2p} \right]                     \\
        \end{array}
    \right.
    \quad \quad , \quad \phi \in \left[ -\frac{\pi}{2},\frac{\pi}{2} \right].
\end{equation}
The sign of the parameter $\phi$ is introduced into equation (\ref{eqn:ParamCrit})$_2$ in order to allow for negative as well as positive values of $\Sigma_m$.
\subsection{Solution of the projection problem for a fixed hardening parameter}
\label{subsec:SolProjProb}
Momentarily assuming the value of the current yield stress $\Sigma$ to be known, we shall now show how the projection problem can be solved through
combination of the yield criterion and the flow rule.
This problem will be reduced to a system of two coupled equations on the unknowns $\phi$ and $\Sigma_{eq}$, which are solved numerically to get
\begin{equation}
\label{eqn:EQ1}
    \fbox{$\displaystyle
      \left[ \frac{6\mu}{3\lambda+2\mu}({\Sigma_m}^* - \Sigma_m) + p\, \Sigma\, \sinh \left( \frac{3}{2}\frac{\Sigma_m}{\Sigma} \right) \right]\Sigma_{eq}
      = p\,\Sigma\, {\Sigma_{eq}}^* \,\sinh \left( \frac{3}{2}\frac{\Sigma_m}{\Sigma} \right).
    $}
    \end{equation}
and 
\begin{equation}
\label{eqn:EQ2}
\fbox{$
  \begin{array}{c}
    \displaystyle \Sigma_{eq}^2 \left\{ 1 + \frac{ A_1{{\mathcal M}_1}^{**} }
        {b^2\left[ \Sigma_{eq} + \frac{3\lambda+2\mu}{45\mu} A_1 \left( {\Sigma_{eq}}^*-\Sigma_{eq} \right) \right]^2 }
        + \frac{ A_2{{\mathcal M}_2}^{**} }
        {b^2\left[ \Sigma_{eq} + \frac{A_2}{5} \left( {\Sigma_{eq}}^*-\Sigma_{eq} \right) \right]^2 }  \right\} \\
    = (1-p)^2 \, {\Sigma}^2 \cos ^2 \phi. \\
  \end{array}
$}
\end{equation}
where ${\Sigma_m}^*$, ${\Sigma_{eq} }^*$, ${ {\mathcal M}_1}^{**}$, and ${ {\mathcal M}_2}^{**}$ are defined as in \textcolor{blue} { \cite {EL09} }. 

{\it This is the second equation of the system looked for} on the unknowns $\phi$ and $\Sigma_{eq}$. The left-hand side depends only on $\Sigma_{eq}$ and the right-hand side only on $\phi$.
$\\$

The simplest way to solve the system of equations \textcolor{blue} { (\ref{eqn:EQ1}, \ref{eqn:EQ2})} on $\phi$ and $\Sigma_{eq}$ may comprise using equation \textcolor{blue} { (\ref{eqn:EQ1})} to express $\Sigma_{eq}$ as a function of $\phi$, and inserting its expression into equation (\ref{eqn:EQ2}) to get an equation on the single unknown $\phi$, to be solved by Newton's method. But numerical experience reveal that the convergence of the Newton iterations is then often problematic. An alternative method comprises solving equation \textcolor{blue} { (\ref{eqn:EQ1})} on $\phi$ through Newton iterations, $\Sigma_{eq}$ being calculated as a function of $\phi$ at each step by solving equation \textcolor{blue}  { (\ref{eqn:EQ2})} through Newton sub-iterations.
\subsubsection{Iterations on the hardening parameter}
\label{subsec:IterHardParam}
The value of the current yield stress $\Sigma$ has been assumed to be known up to now. In reality, it is not and must be determined iteratively. This is done using a fixed point algorithm, starting from the value at time $t$, solving the projection problem with this value, updating it using a discretized form of the evolution equation \textcolor{blue} { (\ref{eqn:EvolEpsBar})}, re-solving the projection problem with the new value, etc. up to convergence.
$\\$

The discretized form of equation \textcolor{blue} {(\ref{eqn:EvolEpsBar})} leads to the following expression of the increment of the hardening parameter $E$:
\begin{equation}
\label{eqn:DeltaE}
  \begin{array}{lll}
    \Delta E & = & \displaystyle \frac{1}{(1-f)\Sigma} \left[ \frac{ \Sigma'_{ij} ({\Sigma'_{ij}}^*-\Sigma'_{ij}) }{2\mu}
                   + \frac{ 3\Sigma_m ({\Sigma_m}^*-\Sigma_m) }{3\lambda+2\mu} \right. \\
    {} & {} & \displaystyle  + \left. \frac{ M'_{ijk} ({M'_{ijk}}^*-M'_{ijk})}{2\mu b^2/5 }
              + \frac{ 3M_{mk} ({M_{mk}}^*-M_{mk}) }{(3\lambda+2\mu)b^2/5 } \right]
  \end{array}
\end{equation}
with, here also, ${\Sigma_m}^*$ and ${M'_{ijk}}^*$ defined as in \textcolor{blue} { \cite {EL09} }. 
\subsubsection{Other features of the numerical implementation}
\label{subsec:OtherFeat}
The numerical implementation of the GLPD model, just like that of all second-gradient models, raises a difficulty tied to the clear necessary use of the second derivatives of the shape functions. This seems to require elements of class ${\mathcal C}^1$ which are never available in standard finite element codes. This difficulty is circumvented through some trick suggested by Gologanu {\it et al.} \textcolor{blue} { \cite{GLPD97}} themselves and used since in several works (see, for instance, Shu {\it et al.} \textcolor{blue} { \cite{SKF99}}, Forest {\it et al.} \textcolor{blue} { \cite{FBC00}}, Matsushima {\it et al.} \textcolor{blue} { \cite{MCC00}} ) for the numerical implementation of various second-gradient models. This trick comprises introducing a new nodal variable in the form of a symmetric second-rank tensor ${\bf W}$, replacing the gradient of the strain rate ${\bf D}$ by the gradient of ${\bf W}$ in all equations and imposing the approximate coincidence of ${\bf W}$ and ${\bf D}$ at the Gauss points through some penalty method.
$\\$

The advantage is that the components of $\nabla{\bf D}$ are then got from the nodal values of the new variable ${\bf W}$ and the sole first derivatives of the shape functions; thus, classical elements of class ${\mathcal C}^0$ are sufficient. The price to pay is an increased number of nodal degrees of freedom: six ($V_1, V_2, W_{11}, W_{22}, W_{12}, W_{33}$) instead of two ($V_1, V_2$) in 2D. Also, imposing the internal constraints $ W_{ij} - D_{ij} = 0 $ by a penalty method may give rise to locking phenomena, for which sub integration is there natural remedy. In practice, 8-node quadratic elements are used with 4-Gauss points integration. Numerical experience reveals that this suffices to prevent locking.
\section{Derivation of the tangent stiffness matrix of the GLPD model}
\label{Ssec: TangentMatrix}
In this section we derive the equations of the tangent stiffness matrix to circumvent global Newton iterations convergence problems the numerical simulations with the GLPD model have revealed.  
We do not compute all the terms of this matrix but only the most important ones. 
Thus, we do not take into account the variations of the stresses due to the variations of the temperature; this is strictly licit because the correction of the stresses is perfomed with explicit scheme, using the stresses at time $t$ and not at time $ t+ \Delta t$; and as a conseqence, the correction is independent upon the displacement increment $ \Delta {\bf U}$ between these instants.
$\\$

Also, we do not take into account the variations of the stresses due to the objective derivation in the law of hypoelasticity, which does indeed depend on $ \Delta {\bf U}$ and therefore generates theoretically a contribution in the stiffness matrix.
Similarly, the influence of geometry on the residual forces will not be taken into account.
We can summarize all this by saying that the calculation of the tangent stifness matrix will be carried out by neglecting the effects of large deformations.
This choice is in conformity with the standard one used in many finite element codes for the calculation of the tangent stiffness matrix for usual elasto-plasticity models (without damage) which practical numerical simulations
involing these models has demonstrated the robustness. 
$\\$

We assume that
\begin{equation}\label{eqn:ContrRed1}
    \left\{
       \begin{array}{lll}
       \tilde \Sigma^\star = \displaystyle \frac{\Sigma^\star}{\bar \Sigma} \quad ; \quad \tilde \Sigma = \displaystyle \frac{\bs}{\bar \Sigma}\\
      {}&{}&{}\\
       \tilde {\bf M}^\star = \displaystyle \frac{{\bf M}^\star}{\bar \Sigma} \quad ; \quad \tilde {\bf M}^{\star \star} = \displaystyle \frac{{\bf M}^{\star \star}}{\bar \Sigma} \quad ; \quad \tilde {\bf M} &=& \displaystyle \frac{{\bf M}}{\bar \Sigma}\\
      {}&{}&{}\\
       \tilde {\bf U} = \displaystyle \bar \Sigma^2 b^2 \frac{\Sigma_{eq}^\star - \Sigma_{eq}}{\Sigma_{eq}} \frac{\bf U}{\bar \Sigma}.\nonumber
       \end{array}
    \right.
\end{equation} 
Hence, we can rewrite the equations giving the expresions $ {\bf U}$, $ \Sigma_{eq}$ and $ \varphi$.
With these notations, we have:
\begin{equation}
\label{eqn:ContrRed2}
    \left\{
       \begin{array}{lll}
       \tilde {M_{ijk}^{\star \star}}' = \displaystyle \tilde {M_{ijk}^{\star}}'- \frac{1}{15} \left(\tilde U_j \delta_{ik} + \tilde U_i \delta_{jk} - \frac{2}{3} \tilde U_k\delta_{ij} \right)\\
      {}&{}&{}\\
       \tilde M_{mk}^{\star\star} = \displaystyle \tilde M_{mk}^{\star} - \frac{3 \lambda + 2 \mu}{45 \mu} \tilde U_k \nonumber
       \end{array}
    \right.
\end{equation} 
and 
\begin{equation} 
\label {eqn:ContrRed3}
   \tilde U_i = \displaystyle \frac{\displaystyle \frac{\tilde {M_{ijj}^{\star}}'}{\tilde \Sigma_{eq} 
          + \displaystyle \frac{A_{II}}{5}({\tilde \Sigma_{eq}^\star-\tilde \Sigma_{eq}})} +\frac{\tilde M_{mi}^{\star}}{\tilde \Sigma_{eq} 
          + \displaystyle\frac{3 \lambda+2 \mu}{45 \mu}A_{I}
             ({\tilde \Sigma_{eq}^\star - \tilde \Sigma_{eq}})}}{ \displaystyle
             \frac{2}{9(\tilde \Sigma_{eq} + \displaystyle \frac{A_{II}}{5}({\tilde \Sigma_{eq}^\star-\tilde \Sigma_{eq}})} + \frac
              {\displaystyle\frac{3 \lambda + 2 \mu}{9 \mu}}{\tilde \Sigma_{eq} + \displaystyle
             \frac{3 \lambda + 2 \mu}{9 \mu}A_{I}(\tilde \Sigma_{eq}^\star-\tilde \Sigma_{eq})} } \equiv \tilde U_i(\tilde \Sigma_{eq}, \tilde \Sigma^\star, \tilde {\bf M}^\star). \nonumber
\end{equation}
The value of $ \tilde U_i $ depends on $ \varphi$ through its arguments $ \tilde \Sigma_{eq}$, $\tilde \Sigma^\star$, and $\tilde {\bf M}^\star$.
Moreover, the equations to be solved on $ \tilde {\Sigma}_{eq}$ and $ \varphi$ are expressed as
\begin{eqnarray}
\label{eqn:ContrRed4}
 \tilde \Sigma_{eq} \left[ 1 + \displaystyle \frac{A_{I}\tilde M_{I}^{\star\star}}{b^2 \left( \tilde \Sigma_{eq} + \displaystyle \frac{3 \lambda + 2 \mu}{45 \mu}A_{I}(\tilde \Sigma_{eq}^\star - \tilde \Sigma_{eq}) \right)^2} + \displaystyle \frac{A_{II}M_{II}^{\star\star}}{ b^2 \left( \tilde \Sigma_{eq} + \displaystyle \frac{A_{II}}{5}({\tilde \Sigma_{eq}^\star-\tilde \Sigma_{eq}}) \right)^2}\right]^\frac{1}{2} \nonumber \\
- (1-p)\cos\varphi = 0, \nonumber
\end{eqnarray}
\begin{equation}
\label{eqn:ContrRed5}
  \frac{6\mu}{3\lambda + 2 \mu}(\tilde \Sigma_m^\star - \tilde \Sigma_m) \tilde {\Sigma}_{eq} - p (\tilde \Sigma_{eq}^\star - \tilde \Sigma_{eq}) \sinh  \left( \frac{3}{2} \tilde \Sigma_m \right) = 0.\nonumber
\end{equation} 
In the subsequent, we introduce the following expressions:
\begin{equation}
\label{eqn:ContrRed6}
   \left\{
       \begin{array}{lll}
        {\mbox {G}} &=& \tilde \Sigma_{eq} \left( 1 + \displaystyle \frac{A_{I}\tilde M_{I}^{\star\star}}{b^2 \left( \tilde \Sigma_{eq} + \displaystyle
               \frac{3\lambda+2 \mu}{45 \mu}A_{I}(\tilde \Sigma_{eq}^\star-\tilde \Sigma_{eq})\right)^2} + 
               \displaystyle \frac{A_{II}\tilde M_{II}^{\star\star}}{b^2 \left( \tilde \Sigma_{eq} + \displaystyle \frac{A_{II}}{5}({\tilde \Sigma_{eq}^\star-\tilde \Sigma_{eq}}) \right)^2} \right)^\frac{1}{2} \\
    &-& (1-p)\cos\varphi \\
       {}&{}&{} \\
       {\mbox {F}} & =& \displaystyle \frac{6\mu}{ 3 \lambda + 2\mu}(\tilde \Sigma_m^\star - \tilde \Sigma_m)\tilde {\Sigma}_{eq} - p(\tilde \Sigma_{eq}^\star - \tilde \Sigma_{eq})\sinh  
\left( \frac{3}{2} \tilde \Sigma_m \right) \nonumber
       \end{array}
    \right.
\end{equation}
where the terms $ {\mbox {G}}$ and $ {\mbox {F}}$ depend on $\varphi$, $\tilde \Sigma_{eq}$, $\tilde {\bf \Sigma}^\star$ , $\tilde {\bf M}^\star$ and $\varphi$, $\tilde {\bf \Sigma}^\star$, $\tilde {\bf M}^\star$.
$\\$

Also, let us assume that 
\begin{equation}
\label{eqn:ContrRed7}
   \left\{
       \begin{array}{lll}
        \mathcal {D}_1 = \tilde \Sigma_{eq} + \displaystyle \frac{3\lambda+2 \mu}{45 \mu}A_{I} \left( \tilde \Sigma_{eq}^\star-\tilde \Sigma_{eq} \right) \quad \equiv \quad \mathcal {D}_1 (\tilde \Sigma_{eq}, \tilde {\bf \Sigma}^\star) \\
       {}&{}&{} \\
        \mathcal {D}_2 = \tilde \Sigma_{eq} + \displaystyle \frac{A_{II}}{45} \left( \tilde \Sigma_{eq}^\star-\tilde \Sigma_{eq} \right) \quad \equiv \quad \mathcal {D}_2(\tilde \Sigma_{eq}, \tilde {\bf \Sigma}^\star) \\
       {}&{}&{} \\
        \mathcal {D} = \displaystyle \frac{2}{9} \frac{1}{\mathcal {D}_2} + \frac{3 \lambda + 2\mu}{45\mu} \frac{1}{\mathcal {D}_1} \quad \equiv \quad \mathcal {D} (\tilde \Sigma_{eq}, \tilde {\bf \Sigma}^\star). \nonumber 
       \end{array}
    \right.
\end{equation}
Then we get 
\begin{equation}
\label{eqn: iContrRed8}
 \tilde U_i = \frac {1}{\mathcal {D}} \left( \frac{\tilde {M_{ijj}^\star}'}{\mathcal {D}_2} + \frac{\tilde M_{mi}^\star}{\mathcal {D}_1} \right). \nonumber
\end{equation}
With Eq.(\ref{eqn: iContrRed8}) we write the expression for $ {\bf G}$ as:
\begin{equation}
\label{eqn:ContrRed6}
{\mbox {G}} = \tilde \Sigma_{eq} \left( 1 + \displaystyle \frac{A_{I} \tilde M_{I}^{\star \star}}{b^2 \mathcal {D}_1^2} + \displaystyle \frac{A_{II}\tilde M_{II}^{\star\star}}{ b^2 \mathcal {D}_2^2} \right)^\frac{1}{2} - ( 1 - p ) \cos \varphi. \nonumber
\end{equation}
\subsection{Calculation of the derivatives}
\label{deriv}
First, for practical purposes, we find the following intermediate expressions:
\begin{eqnarray}
\label{eqn:Deriv1}
    \left\{
       \begin{array}{lll}
    \displaystyle \frac{\partial \tilde \Sigma_{eq}^\star}{\partial \tilde \Sigma_{ij}^\star} &=& \displaystyle \frac{3}{2}\frac{\tilde {\Sigma_{ij}^\star}'}{\tilde \Sigma_{eq}^\star} \nonumber\\
    \tilde M_{mk}^\star &=& \displaystyle \frac{1}{3} \tilde M_{hhk}^\star \Rightarrow \frac{\partial \tilde M_{mk}^\star}{\partial \tilde M_{pqr}} = \frac{1}{3} \delta_{pq}\delta_{kr} \nonumber\\
       {}&{}&{} \\
    \tilde {M_{ijj}^\star}' &=& \displaystyle \tilde M_{ijj}^\star - \frac{1}{3}\delta_{ij}\tilde M_{hhj}^\star = \tilde M_{ijj}^\star - \frac{1}{3}\tilde M_{hhi}^\star \nonumber\\
       {}&{}&{} \\
    \Rightarrow \displaystyle \frac{\partial \tilde {M_{ijj}^\star}'}{\partial \tilde {M_{pqr}^\star}'} &=& \displaystyle \frac{1}{2}\left( \delta_{ip}\delta_{jq} + \delta_{iq}\delta_{jp}\right)\delta_{jr} - \frac{1}{3}\delta_{pq}\delta_{ir} \nonumber\\
       {}&{}&{} \\
    &=& \displaystyle \frac{1}{2}\delta_{ip}\delta_{qr} + \frac{1}{2}\delta_{iq}\delta_{pr} - \frac{1}{3}\delta_{pq}\delta_{ir}, \nonumber
       \end{array}
    \right.
\end{eqnarray}
$\\$

\begin{equation}
\label{eqn:Deriv2}
    \left\{
       \begin{array}{lll}
        \displaystyle \frac{\partial \mathcal {D}_1}{\partial \tilde \Sigma_{eq}} = 1 - \displaystyle \frac{3\lambda + 2 \mu}{45 \mu}A_{I} \quad;\quad 
        \displaystyle \frac{\partial \mathcal {D}_2}{\partial \tilde \Sigma_{eq}} = 1 - \displaystyle \frac{A_{II}}{5}\\
       {}&{}&{} \\
        \displaystyle \frac{\partial \mathcal D_1}{\partial \tilde \Sigma_{pq}^\star} = \displaystyle \frac{3\lambda + 2 \mu}{45 \mu}A_{I}\frac{3}{2}\frac{\tilde {\Sigma_{pq}^\star}'}{\tilde \Sigma_{eq}^\star} \quad;\quad \displaystyle \frac{\partial \mathcal {D}_2}{\partial \tilde \Sigma_{pq}^\star} = \displaystyle \frac{A_{II}}{5}\frac{3}{2}\frac{\tilde {\Sigma_{pq}^\star}'}{\tilde \Sigma_{eq}^\star} \\
       {}&{}&{} \\
        \displaystyle \frac{\partial \mathcal {D} }{\partial \tilde \Sigma_{eq}} = -\displaystyle \frac{3\lambda + 2 \mu}{45 \mu} \frac{1}{\mathcal {D}_1^2}\frac{\partial \mathcal {D}_1}{\partial \tilde \Sigma_{eq}} - \frac{2}{9}\frac{1}{\mathcal {D}_2^2}\frac{\partial \mathcal {D}_2}{\partial \tilde \Sigma_{eq}}\\
       {}&{}&{} \\
        \displaystyle \frac{\partial \mathcal {D}}{\partial \tilde \Sigma_{pq}^\star} = -\displaystyle \frac{3\lambda + 2 \mu}{45 \mu} \frac{1}{\mathcal {D}_1^2}\frac{\partial \mathcal {D}_1}{\partial \tilde \Sigma_{pq}} - \frac{2}{9}\frac{1}{\mathcal {D}_2^2}\frac{\partial \mathcal {D}_2}{\partial \tilde \Sigma_{pq}}.\nonumber
       \end{array}
    \right.
\end{equation}
\subsubsection{Derivatives of the $ \tilde U_i$}
\label{deriv}
The terms $ \tilde  U_i$ depends on $\tilde \Sigma_{eq}$, $\tilde {\bf \Sigma}^\star$ et $\tilde {\bf M}^\star$; as a consequence, we have
\begin{eqnarray}
\label{eqn:Deriv1}
    \frac{\partial \tilde U_i}{\partial \tilde \Sigma_{eq}} &=& - \frac{1}{\mathcal {D}^2}\frac{\partial \mathcal {D}}{\partial \tilde \Sigma_{eq}} \left( \frac{\tilde {M_{ijj}^\star}'}{\mathcal {D}_2} + \frac{\tilde M_{mi}^\star}{\mathcal {D}_1}\right) - \frac{1}{\mathcal {D}}  \left( \frac{ \tilde {M_{ijj}^\star}'}{\mathcal {D}_2^2}\frac{\partial \mathcal {D}_2}{\partial \tilde \Sigma_{eq}} + \frac{\tilde M_{mi}^\star}{\mathcal {D}_1^2}\frac{\partial \mathcal {D}_1}{\partial \tilde \Sigma_{eq}}\right)\nonumber\\
    &=& -\left( \frac{1}{\mathcal {D}^2 \mathcal {D}_2} \frac{\partial \mathcal {D}}{\partial \tilde \Sigma_{eq}} + \frac{1}{\mathcal {D} \mathcal {D}_2^2} \frac{\partial \mathcal {D}_2}{\partial \tilde \Sigma_{eq}}\right)\tilde {M_{ijj}^\star}' - \left( \frac{1}{\mathcal {D}^2 \mathcal {D}_1} \frac{\partial \mathcal {D}}{\partial \tilde \Sigma_{eq}} + \frac{1}{\mathcal {D} \mathcal {D}_1^2} \frac{\partial \mathcal {D}_1}{\partial \tilde \Sigma_{eq}}\right)\tilde M_{mi}^\star \quad ; \quad \nonumber\\
    {}&{}&{} \nonumber\\ 
\frac{\partial \tilde U_i}{\partial \tilde \Sigma_{pq}^\star} &=& -\left( \frac{1}{\mathcal {D}^2 \mathcal {D}_2} \frac{\partial \mathcal {D}}{\partial \tilde \Sigma_{pq}^\star} + \frac{1}{\mathcal {D} \mathcal {D}_2^2} \frac{\partial \mathcal {D}_2}{\partial \tilde \Sigma_{pq}^\star} \right)\tilde {M_{ijj}^\star}' - \left( \frac{1}{\mathcal {D}^2 \mathcal {D}_1} \frac{\partial \mathcal {D}}{\Sigma_{pq}} + \frac{1}{\mathcal {D} \mathcal {D}_1^2} \frac{\partial \mathcal {D}_1}{\partial \tilde \Sigma_{pq}^\star}\right)\tilde M_{mi}^\star  \quad \mbox {and} \quad \nonumber\\
    {}&{}&{} \nonumber\\
\frac{\partial \tilde U_i}{\partial \tilde M_{pqr}^\star} &=& \frac{1}{\mathcal {D}} \left( \frac{1}{\mathcal {D}_1} \frac{\partial \tilde M_{mi}^\star}{\partial \tilde M_{pqr}^\star} + \frac{1}{\mathcal {D}_2} \frac{\partial \tilde M_{ijj}^\star}{\partial \tilde M_{pqr}^\star}\right) \nonumber\\
    &=& \frac{1}{\mathcal {D}} \left[ \frac{1}{\mathcal {D}_1}\frac{1}{3} \delta_{pq}\delta_{ir} + \frac{1}{\mathcal D_2}(\frac{1}{2}\delta_{ip}\delta_{qr} + \frac{1}{2}\delta_{iq}\delta_{pr} - \frac{1}{3}\delta_{pq}\delta_{ir}) \right].\nonumber
\end{eqnarray}
\subsubsection{The derivatives of the terms $ \tilde {M_{ijk}^{\star \star}}'$ and $ \tilde M_{mk}^{\star \star}$}
\label{deriv2}
The terms $ \tilde {M_{ijk}^{\star \star}}'$ and $ \tilde M_{mk}^{\star \star}$ also depend on $(\tilde \Sigma_{eq}, \tilde {\bf \Sigma}^\star, \tilde {\bf M}^\star)$; thus we get 
\begin{equation}
\label{eqn:Deriv3}
    \left\{
       \begin{array}{lll}
       \displaystyle \frac{\partial \tilde {M_{ijk}^{\star \star}}'}{\partial \tilde \Sigma_{eq}} &=& - \displaystyle \frac{1}{15} \left(\frac{\partial \tilde U_j}{\partial \tilde \Sigma_{eq}} \delta_{ik} + \frac{\partial \tilde U_i}{\partial \tilde \Sigma_{eq}}\delta_{jk} - \frac{2}{3} \frac{\partial \tilde U_k}{\partial \tilde \Sigma_{eq}}\delta_{ij} \right)\\
       {}&{}&{} \\
       \displaystyle \frac{\partial \tilde M_{mk}^{\star\star}}{\partial \tilde \Sigma_{eq}} &=&  - \displaystyle \frac{3\lambda+2\mu}{45\mu}\frac{\partial \tilde U_k}{\partial \tilde \Sigma_{eq}}, \nonumber
       \end{array}
    \right.
\end{equation}
\begin{equation}
\label{eqn:Deriv4}
    \left\{
       \begin{array}{lll}
       \displaystyle \frac{\partial \tilde {M_{ijk}^{\star \star}}'}{\partial \tilde \Sigma_{pq}^\star} &=& - \displaystyle \frac{1}{15} \left(\frac{\partial \tilde U_j}{\partial \tilde \Sigma_{pq}^\star} \delta_{ik} + \frac{\partial \tilde U_i}{\partial \tilde \Sigma_{pq}^\star}\delta_{jk} - \frac{2}{3} \frac{\partial \tilde U_k}{\partial \tilde \Sigma_{pq}^\star}\delta_{ij} \right)\\
       {}&{}&{} \\
       \displaystyle \frac{\partial \tilde M_{mk}^{\star\star}}{\partial \tilde \Sigma_{pq}^\star} &=&  - \displaystyle \frac{3\lambda+2\mu}{45\mu}\frac{\partial \tilde U_k}{\partial \tilde \Sigma_{pq}^\star}\nonumber
       \end{array}
    \right.
\end{equation}
and 
\begin{equation}
\label{eqn:Deriv5}
    \left\{
       \begin{array}{lll}
       \displaystyle \frac{\partial \tilde {M_{ijk}^{\star \star}}'}{\partial \tilde M_{pqr}^\star} &=& \displaystyle \frac{1}{2} (\delta_{ip}\delta_{jq} + \delta_{iq}\delta_{jp})\delta_{kr} - \frac{1}{3}\delta_{ij}\delta_{pq}\delta_{kr} \\ 
                    &-& \displaystyle \frac{1}{15} \left(\frac{\partial \tilde U_j}{\partial \tilde M_{pqr}^\star} \delta_{ik} + \frac{\partial \tilde U_i}{\partial \tilde M_{pqr}^\star}\delta_{jk} - \frac{2}{3} \frac{\partial \tilde U_k}{\partial \tilde M_{pqr}^\star}\delta_{ij} \right)\\
       {}&{}&{} \\
       \displaystyle \frac{\partial \tilde M_{mk}^{\star\star}}{\partial \tilde M_{pqr}^\star} &=& \displaystyle \frac{1}{3}\delta_{pq}\delta_{kr} - \displaystyle \frac{3\lambda+2\mu}{45\mu}\frac{\partial \tilde U_k}{\partial \tilde M_{pqr}^\star}.\nonumber
       \end{array}
    \right.
\end{equation}
\subsubsection{The derivatives of $ \tilde M_{I}^{\star \star}$ and $ \tilde M_{II}^{\star \star}$}
\label{deriv}
The terms $ \tilde M_{I}^{\star \star}$ and $ \tilde M_{II}^{\star \star}$ depend on the variables $(\tilde \Sigma_{eq}, \tilde {\bf \Sigma}^ \star, \tilde {\bf M}^\star)$.
Taking the derivatives we get
\begin{equation}
\label{eqn:Deriv6}
    \left\{
       \begin{array}{lll}
       \tilde M_{I}^{\star \star} = \displaystyle \frac{1}{9}\tilde M_{hhi}^{\star\star}\tilde M_{kki}^{\star\star} \quad \Rightarrow \quad \displaystyle \frac{\partial \tilde M_{I}^{\star \star}}{\partial \tilde \Sigma_{eq}} = \frac{2}{9} \tilde M_{kki}^{\star \star} \frac{\partial \tilde M_{hhi}^{\star \star}}{\partial \tilde \Sigma_{eq}} = \displaystyle 2 \tilde M_{mk}^{\star \star} \frac{\partial \tilde M_{mk}^{\star \star}}{\partial \tilde \Sigma_{eq}}\\
       {}&{}&{} \\
       \tilde M_{II}^{\star \star} = \displaystyle \frac{3}{2}\tilde {M_{ijk}^{\star \star}}'\tilde {M_{ijk}^{\star \star}}' \quad \Rightarrow \quad \displaystyle \frac{\partial \tilde M_{II}^{\star \star}}{\partial \tilde \Sigma_{eq}} = 3 \tilde {M_{ijk}^{\star \star}}' \frac{\partial \tilde {M_{ijk}^{\star \star}}'}{\partial \tilde \Sigma_{eq}} = \displaystyle -\frac{2}{5} \tilde {M_{ijj}^{\star \star}}' \frac{\partial \tilde U_i}{\partial \tilde \Sigma_{eq}},\nonumber
       \end{array}
    \right.
\end{equation}
\begin{equation}
\label{eqn:Deriv6}
    \left\{
       \begin{array}{lll}
       \displaystyle \frac{\partial \tilde M_{I}^{\star \star}}{\partial \tilde \Sigma_{pq}^\star} = \frac{2}{9} \tilde M_{hhi}^{\star \star} \frac{\partial \tilde M_{kki}^{\star \star}}{\partial \tilde \Sigma_{pq}^\star} = \displaystyle 2 \tilde M_{mk}^{\star \star} \frac{\partial \tilde M_{mk}^{\star \star}}{\partial \tilde \Sigma_{pq}^\star}\\
       {}&{}&{} \\
       \displaystyle \frac{\partial \tilde M_{II}^{\star \star}}{\partial \tilde \Sigma_{pq}^\star} = 3 \tilde {{M_{ijk}^{\star \star}}'} \frac{\partial \tilde {M_{ijk}^{\star \star}}'}{\partial \tilde \Sigma_{pq}^\star} = \displaystyle -\frac{2}{5} \tilde {M_{ijj}^{\star \star}}' \frac{\partial \tilde U_i}{\partial \tilde \Sigma_{pq}^\star}, \nonumber\\
       \end{array}
     \right.
\end{equation}
and 
\begin{equation}
\label{eqn:Deriv6}
    \left\{
       \begin{array}{lll}
       \displaystyle \frac{\partial \tilde M_{I}^{\star \star}}{\partial \tilde M_{pqr}^\star} = \frac{2}{9} \tilde M_{hhi}^{\star \star} \frac{\partial \tilde M_{kki}^{\star \star}}{\partial \tilde M_{pqr}^\star} = \displaystyle 2 \tilde M_{mk}^{\star \star} \frac{\partial \tilde M_{mk}^{\star \star}}{\partial \tilde M_{pqr}^\star}\\
       {}&{}&{} \\
       \displaystyle \frac{\partial \tilde M_{II}^{\star \star}}{\partial \tilde M_{pqr}^\star} = 3 \tilde {M_{ijk}^{\star \star}}' \frac{\partial \tilde {M_{ijk}^{\star \star}}'}{\partial \tilde M_{pqr}^\star} = 3 \tilde {M_{pqr}^{\star \star}}' -\displaystyle \frac{2}{5} \tilde {M_{ijj}^{\star \star}}' \frac{\partial \tilde U_i}{\partial \tilde M_{pqr}^\star}. \nonumber
       \end{array}
    \right.
\end{equation}
\subsubsection{The derivatives of $ {\mbox {G}}$ }
\label{deriv}
We know that $ {\mbox {G}} \equiv {\bf G} (\varphi, \tilde \Sigma_{eq}, \tilde {\bf \Sigma}^\star, \tilde {\bf M}^\star)$. by posing that 
\begin{equation}
\label{eqn:DerivG1}
\mathcal {S} =  \left( 1 + \displaystyle \frac{A_{I}\tilde M_{I}^{\star \star}}{b^2 \mathcal D_1^2} + 
\displaystyle \frac{A_{II}\tilde M_{II}^{\star\star}}{ b^2 \mathcal D_2^2} \right)^\frac{1}{2} \equiv \mathcal {S}(\tilde \Sigma_{eq}, \tilde {\bf \Sigma}^\star, \tilde {\bf M}^\star), \nonumber
\end{equation}
we get
\begin{equation}
\label{eqn:1DerivG1}
{\mbox {G}} = \tilde \Sigma_{eq}\mathcal {S} - (1-p) \cos \varphi. \nonumber
\end{equation}
Therafter
\begin{equation}
\label{eqn:DerivG2}
    \left\{
       \begin{array}{lll}
       \displaystyle \frac{\partial {\mbox {G}}}{\partial \tilde \Sigma_{eq}} &=& \mathcal {S} + \displaystyle \tilde \Sigma_{eq}\frac{\partial \mathcal {S}}{\partial \tilde \Sigma_{eq}} \\
           {}&{}&{}\\
&=&  \displaystyle \mathcal {S} + \displaystyle \frac{\tilde \Sigma_{eq}}{2\mathcal {S}}\left( \frac{A_{I}}{b^2 \mathcal {D}_1^2}\frac{\partial \tilde M_{I}^{\star \star}}{\partial \tilde \Sigma_{eq}} + \frac{A_{II}}{b^2 \mathcal {D}_2^2}\frac{\partial \tilde M_{II}^{\star \star}}{\partial \tilde \Sigma_{eq}} - \frac{2A_{II}}{b^2 \mathcal {D}_1^3}\tilde M_{II}^{\star \star} \frac{\partial \mathcal {D}_1}{\partial \tilde \Sigma_{eq}} - \frac{2A_{II}}{b^2\mathcal {D}_2^3}\tilde M_{II}^{\star \star} \frac{\partial \mathcal {D}_2}{\partial \tilde \Sigma_{eq}}\right)\\
           {}&{}&{}\\
       \displaystyle \frac{\partial {\mbox {G}}}{\partial \tilde \Sigma_{pq}^\star} 
&=&  \displaystyle \tilde \Sigma_{eq} \frac{\partial \mathcal {S}}{\partial \tilde \Sigma_{pq}^\star}\\
          {}&{}&{}\\
 &=&  \displaystyle \displaystyle \frac{\tilde \Sigma_{eq}}{2\mathcal {S}}\left( \frac{A_{I}}{b^2 \mathcal {D}_1^2}\frac{\partial \tilde M_{I}^{\star \star}}{\partial \tilde \Sigma_{pq}^\star} + \frac{A_{II}}{b^2\mathcal {D}_2^2}\frac{\partial \tilde M_{II}^{\star \star}}{\partial \tilde \Sigma_{pq}^\star} - \frac{2A_{I}}{b^2\mathcal {D}_1^3}\tilde M_{I}^{\star \star} \frac{\partial \mathcal {D}_1}{\partial \tilde \Sigma_{pq}^\star} - \frac{2A_{II}}{b^2\mathcal {D}_2^3}\tilde M_{II}^{\star \star} \frac{\partial \mathcal {D}_2}{\partial \tilde \Sigma_{pq}^\star}\right)\\
           {}&{}&{}\\
       \displaystyle \frac{\partial {\mbox {G}}}{\partial \tilde M_{pqr}^\star} &=&  \displaystyle \tilde \Sigma_{eq} \frac{\partial \mathcal {S}}{\partial \tilde M_{pqr}^\star} =  \displaystyle \displaystyle \frac{\tilde \Sigma_{eq}}{2\mathcal {S}}\left( \frac{A_{I}}{b^2 \mathcal {D}_2^2}\frac{\partial \tilde M_{I}^{\star \star}}{\partial \tilde M_{pqr}^\star} + \frac{A_{II}}{ b^2 \mathcal {D}_2^2}\frac{\partial \tilde M_{II}^{\star \star}}{\partial \tilde M_{pqr}^\star} \right)\\
           {}&{}&{}\\
       \displaystyle \frac{{\partial \mbox {G}}}{\partial \varphi} &=& ( 1-p )\sin \varphi.\nonumber
       \end{array}
    \right.
\end{equation}
With the derivatives of $ {\mbox {G}}$ we can find the derivatives of $ \tilde \Sigma_{eq}$ with respect to $ \varphi$, $ \tilde \Sigma^\star$, $ \tilde {\bf M}^\star$.
\subsubsection{The derivatives of the term $ \tilde \Sigma_{eq}$ with respect to $ \varphi$, $ \tilde \Sigma^\star$, $ \tilde {\bf M}^\star$}
\label{deriv}
The following implications hold:
\begin{eqnarray}
\label{eqn:DerivG3}
    {\mbox {G}} \equiv {\mbox {G}}(\varphi, \tilde \Sigma_{eq}, \tilde {\bf \Sigma}^\star, \tilde {\bf M}^\star) = 0 \quad &\Rightarrow& \quad \frac{\partial {\mbox {G}}}{\partial \tilde \Sigma_{eq}} D \tilde \Sigma_{eq} + \frac{{\partial \mbox {G}}}{\partial \varphi} D \varphi + \frac{\partial {\bf G}}{\partial \tilde \Sigma_{pq}^\star} D \tilde \Sigma_{pq}^\star + \frac{\partial {\mbox {G}}}{\partial \tilde M_{pqr}^\star} D \tilde M_{pqr}^\star = 0 \nonumber\\
   &\Rightarrow& \frac{\partial \tilde \Sigma_{eq}}{\partial \varphi} = - \frac{\partial {\mbox {G}}/\partial \varphi}{\partial {\mbox {G}}/\partial \tilde \Sigma_{eq}} \quad ; \quad \frac{\partial \tilde \Sigma_{eq}}{\partial \tilde \Sigma_{pq}^\star} = -\frac{\partial {\mbox {G}}/\partial \tilde \Sigma_{pq}^\star}{\partial {\mbox {G}}/\partial \tilde \Sigma_{eq}}; \nonumber\\
\quad &\quad& \quad \frac{\partial \tilde \Sigma_{eq}}{\partial \tilde M_{pqr}^\star} = -\frac{\partial {\mbox {G}}/\partial \tilde M_{pqr}^\star}{\partial {\mbox {G}}/\partial \tilde \Sigma_{eq}}.\nonumber
\end{eqnarray}
%%
%%%%%%%%%%%%%%%%%%%%%%%%%
\subsubsection{The derivative of $ {\mbox {F}}( \varphi, \tilde \Sigma^\star,  \tilde {\bf M}^\star) $}
\label{sec: deriv11d}
From the formula defining $ {\mbox {F}}$, we immediately obtain the derivatives
\begin{equation}
\label{eqn:DerivG4}
    \left\{
       \begin{array}{lll}
       \displaystyle \frac{\partial {\mbox {F}}}{\partial \varphi} &=& \displaystyle \left[ - \displaystyle \frac{6\mu}{3\lambda + 2\mu}\tilde \Sigma_{eq}  - \displaystyle \frac{3}{2}p(\tilde \Sigma_{eq}^\star - \tilde \Sigma_{eq}) \cosh ( \displaystyle \frac{3}{2} \tilde \Sigma_m) \right] \displaystyle \frac{D \tilde \Sigma_m}{D \varphi}  \\
           {}&{}&{}\\
& + & \displaystyle \left[ \frac{6\mu}{3\lambda + 2\mu}(\tilde \Sigma_m^\star - \tilde \Sigma_m) + p \sinh ( \frac{3}{2} \tilde \Sigma_m ) \right]\frac{\partial \tilde \Sigma_{eq}}{\partial \varphi}\\
           {}&{}&{}\\
       \displaystyle \frac{\partial {\mbox {F}}}{\partial \tilde \Sigma_{pq}^\star} &=& \displaystyle \frac{2\mu}{3\lambda + 2\mu} \delta_{pq} \tilde \Sigma_{eq} - \frac{3p}{2} \sinh (\frac{3}{2} \tilde \Sigma_m) \frac{ \tilde {\Sigma_{pq}^\star}'}{\tilde \Sigma_{eq}^\star} \\
           {}&{}&{}\\
      & + &  \displaystyle \left[ \frac{6\mu}{3\lambda + 2\mu}(\tilde \Sigma_m^\star - \tilde \Sigma_m) + p \sinh ( \frac{3}{2} \tilde \Sigma_m ) \right]\frac{\partial \tilde \Sigma_{eq}}{\partial \tilde \Sigma_{pq}^\star}\\
           {}&{}&{}\\
       \displaystyle \frac{\partial {\mbox {F}}}{\partial \tilde M_{pqr}^\star} &=&  \displaystyle \left[ \frac{6\mu}{3\lambda + 2\mu}(\tilde \Sigma_m^\star - \tilde \Sigma_m) + p \sinh ( \frac{3}{2} \tilde \Sigma_m ) \right]\frac{\partial \tilde \Sigma_{eq}}{\partial \tilde M_{pqr}^\star}\nonumber
       \end{array} 
    \right.
\end{equation}
which enable thre calculation of the derivatives of the term $ \varphi$ with respect to $\tilde \Sigma^\star$ and $\tilde {\bf M}^\star$.
\subsubsection{The derivatives of the terms $ \varphi(\tilde \Sigma^\star, \tilde {\bf M}^\star)$}
\label{deriv20}
Following the results obtained in Section ( \ref{sec: deriv11d}), we have  
\begin{eqnarray}
\label{eqn:DerivG5}
    {\mbox {F}} \equiv {\mbox {F}}(\varphi, \tilde {\bf \Sigma}^\star, \tilde {\bf M}^\star) = 0 \quad &\Rightarrow& \quad \frac{\partial {\mbox {F}}}{\partial \varphi } D \varphi  + \frac{{\partial \mbox {F}}}{\partial \tilde \Sigma_{pq}^\star} D \tilde \Sigma_{pq}^\star  + \frac{\partial {\mbox {F}}}{\partial \tilde M_{pqr}^\star} D \tilde M_{pqr}^\star = 0 \nonumber\\
   &\Rightarrow& \frac{\partial \varphi }{\partial \tilde \Sigma_{pq}^\star} = -\frac{\partial {\mbox {F}}/\partial \tilde \Sigma_{pq}^\star}{\partial {\mbox {F}}/\partial \varphi} \quad ; \quad \frac{\partial \varphi}{\partial \tilde M_{pqr}^\star} = -\frac{\partial {\mbox {F}}/\partial \tilde M_{pqr}^\star}{\partial {\mbox {F}}/\partial \varphi}. \nonumber
\end{eqnarray}
\subsubsection{The derivatives of $ \tilde \Sigma_{eq}(\tilde \Sigma^\star, \tilde {\bf M}^\star)$}
\label{deriv21}
We denote the derivatives in this section by $ \displaystyle \frac{D \tilde \Sigma_{eq}}{D \tilde \Sigma_{pq}^\star}$ and $ \displaystyle \frac{\mbox { D} \tilde \Sigma_{eq}}{D \tilde M_{pqr}^\star}$.
We get
\begin{equation}
\label{eqn:DerivG5}
    \left\{
       \begin{array}{lll}
         \displaystyle \frac{D \tilde \Sigma_{eq}}{D \tilde \Sigma_{pq}^\star} = \displaystyle \frac{\partial \tilde \Sigma_{eq}}{\partial \tilde \Sigma_{pq}^\star} + \frac{\partial \tilde \Sigma_{eq}}{\partial \varphi} \frac{\partial \varphi}{\partial \tilde \Sigma_{pq}^\star}\\
           {}&{}&{}\\
         \displaystyle \frac{D \tilde \Sigma_{eq}}{D \tilde M_{pqr}^\star} = \displaystyle \frac{\partial \tilde \Sigma_{eq}}{\partial \tilde M_{pqr}^\star} + \frac{\partial \tilde \Sigma_{eq}}{\partial \varphi} \frac{\partial \varphi}{\partial \tilde M_{pqr}^\star}.\nonumber
       \end{array}
    \right.
\end{equation}
\subsubsection{The derivatives of $ \tilde \Sigma_m (\tilde \Sigma^\star, \tilde {\bf M}^\star )$}
\label{deriv22}
We obtain 
\begin{eqnarray}
\label{eqn:DerivG6}
\frac{\partial \tilde \Sigma_m}{\partial \tilde \Sigma_{pq}^\star} = \frac{D \tilde \Sigma_m}{D \varphi} \frac{\partial \varphi}{\partial \Sigma_{pq}^\star} \quad ; \quad \frac{\partial \tilde \Sigma_m}{\partial \tilde M_{pqr}^\star} = \frac{D \tilde \Sigma_m}{D \varphi} \frac{\partial \varphi}{\partial M_{pqr}^\star}.\nonumber
\end{eqnarray}
\subsubsection{The derivatives of $ \tilde {\bs}' (\tilde \Sigma^\star; \tilde {\bf M}^\star)$}
\label{deriv1}
%%%%%%%%%%%%%%%%%%%%%%%%%
%%
Also, we obtain
\begin{eqnarray}\label{eqn:DerivG7}
\Sigma'_{ij} = \frac{\tilde \Sigma_{eq}}{\tilde \Sigma_{eq}^\star} {\Sigma_{ij}^\star}' \quad \Rightarrow \quad \frac{\partial \tilde \Sigma'_{ij}}{\partial \tilde \Sigma_{pq}^\star} = \frac{1}{\tilde \Sigma_{eq}^\star}\frac{D \tilde \Sigma_{eq}}{D \tilde \Sigma_{pq}^\star}\tilde {\Sigma_{ij}^\star}' - \frac{\tilde \Sigma_{eq}}{(\tilde \Sigma_{eq}^\star)^2} \frac{\partial \tilde \Sigma_{eq}^\star}{\partial \tilde \Sigma_{pq}^\star} \tilde {\Sigma_{ij}^\star}' + \frac{\tilde \Sigma_{eq}}{\tilde \Sigma_{eq}^\star}\frac{\partial \tilde {\Sigma_{ij}^\star}'}{\tilde \Sigma_{pq}^\star} \nonumber
\end{eqnarray}
where
\begin{eqnarray}
\label{eqn:DerivG8}
    \frac{\partial \tilde \Sigma_{eq}^\star}{\partial \tilde \Sigma_{pq}^\star} &=& \frac{3}{2} \frac{\tilde {\Sigma_{pq}^\star}'}{\tilde \Sigma_{eq}^\star} \quad \mbox {et} \quad \frac{\partial \tilde {\Sigma_{ij}^\star}'}{\partial \tilde \Sigma_{pq}^\star} = \frac{1}{2}(\delta_{ip}\delta_{jq} + \delta_{iq}\delta_{jp}) - \frac{1}{3}\delta_{ij}\delta_{pq} \nonumber\\
  \frac{\partial \tilde \Sigma'_{ij}}{\partial \tilde M_{pqr}^\star} &=& \frac{1}{\tilde \Sigma_{eq}^\star}\frac{D \tilde \Sigma_{eq}}{D \tilde M_{pqr}^\star}\tilde {\Sigma_{ij}^\star}'. \nonumber
\end{eqnarray}
It becomes possible, from these last equations, to easily compute $ \displaystyle \frac{\partial \tilde \Sigma_{ij}}{\partial \tilde \Sigma_{pq}^\star}$ and $ \displaystyle \frac{\partial \tilde \Sigma_{ij}}{\partial \tilde M_{pqr}^\star}$.
\subsubsection{The derivatives of the terms $ \tilde M_{mk}^{\star \star} (\tilde \Sigma^\star; \tilde {\bf M}^\star)$ et $ \tilde M_{ijk}'^{\star \star} (\tilde \Sigma^\star; \tilde {\bf M}^\star)$}
\label{deriv24}
The derivatives of $ \tilde M_{mk}^{\star \star}$ is given by
\begin{equation}
\label{eqn:DerivG5}
    \left\{
       \begin{array}{lll}
         \displaystyle \frac{D \tilde M_{mk}^{\star \star}}{D \tilde \Sigma_{pq}^\star} = \displaystyle \frac{\partial \tilde M_{mk}^{\star \star}}{\partial \tilde \Sigma_{pq}^\star} + \frac{\partial \tilde M_{mk}^{\star \star}}{\partial \tilde \Sigma_{eq}}\frac{D \tilde \Sigma_{eq}}{D \tilde \Sigma_{pq}^\star} \quad ; \quad 
         \displaystyle \frac{D \tilde M_{mk}^{\star \star}}{D \tilde M_{pqr}^\star} = \displaystyle \frac{\partial \tilde M_{mk}^{\star \star}}{\partial \tilde M_{pqr}^\star} + \frac{\partial \tilde M_{mk}^{\star \star}}{\partial \tilde \Sigma_{eq}}\frac{D \tilde \Sigma_{eq}}{D \tilde M_{pqr}^\star}\\
           {}&{}&{}\\
         \displaystyle \frac{D \tilde {M_{ijk}^{\star \star}}'}{D \tilde \Sigma_{pq}^\star} = \displaystyle \frac{\partial \tilde {M_{ijk}^{\star \star}}'}{\partial \tilde \Sigma_{pq}^\star} + \frac{\partial \tilde {M_{ijk}^{\star \star}}'}{\partial \tilde \Sigma_{eq}}\frac{D \tilde \Sigma_{eq}}{D \tilde \Sigma_{pq}^\star} \quad ; \quad 
         \displaystyle \frac{D \tilde {M_{ijk}^{\star \star}}'}{D \tilde M_{pqr}^\star} = \displaystyle \frac{\partial \tilde {M_{ijk}^{\star \star}}'}{\partial \tilde M_{pqr}^\star} + \frac{\partial \tilde {M_{ijk}^{\star \star}}'}{\partial \tilde \Sigma_{eq}}\frac{D \tilde \Sigma_{eq}}{D \tilde M_{pqr}^\star}.\nonumber
       \end{array}
    \right.
\end{equation}
\subsubsection{The derivative of the term $ \tilde M_{mk} \equiv \tilde M_{mk} (\tilde \Sigma^\star; \tilde {\bf M}^\star)$ }
\label{deriv}
%%%%%%%%%%%%%%%%%%%%%%%%%
%%
We know that
\begin{equation}
\label{eqn:ParamCrit2}
\tilde M_{mk} = \displaystyle R_1\tilde M_{mk}^{\star \star} \quad \mbox {where} \quad R_1 = \displaystyle \sqrt{\frac{M_{I}}{M_{I}^{ \star \star}}} = \frac{\tilde \Sigma_{eq}}{\tilde \Sigma_{eq} + \displaystyle \frac{3\lambda+2\mu}{45 \mu}A_{I} (\tilde \Sigma_{eq}^\star - \tilde \Sigma_{eq})} = \frac{\tilde \Sigma_{eq}}{\mathcal {D}_1}; \nonumber
\end{equation}
as a consequence we get  
\begin{equation}\label{eqn:ParamCrit2}
    \left\{
       \begin{array}{lll}
           \displaystyle \frac{\partial R_1}{\partial \tilde \Sigma_{pq}^\star} = \displaystyle \frac{1}{\mathcal D_1}\frac{D \tilde \Sigma_{eq}}{D \tilde \Sigma_{pq}^\star} - \frac{\tilde \Sigma_{eq}}{\mathcal D_1^2}\left( \frac{\partial \mathcal {D}_1}{\partial \tilde \Sigma_{pq}^\star} + \frac{\partial \mathcal {D}_1}{\partial \tilde \Sigma_{eq}}\frac{D \tilde \Sigma_{eq}}{D \tilde \Sigma_{pq}^\star}\right)\\
           {}&{}&{}\\
           \displaystyle \frac{\partial R_1}{\partial \tilde M_{pqr}^\star} = \displaystyle \frac{1}{\mathcal {D}_1}\frac{D \tilde \Sigma_{eq}}{D \tilde M_{pqr}^\star} - \frac{\tilde \Sigma_{eq}}{\mathcal {D}_1^2} \frac{\partial \mathcal {D}_1}{\partial \tilde \Sigma_{eq}}\frac{D \tilde \Sigma_{eq}}{D \tilde M_{pqr}^\star}.\nonumber 
       \end{array}
    \right.
\end{equation}
Thereafter, 
\begin{equation}\label{eqn:ParamCrit2}
     \displaystyle \frac{\partial \tilde M_{mk}}{\partial \Sigma_{pq}^\star} = R_2 \displaystyle \frac{D \tilde M_{mk}^{\star \star}}{D \tilde \Sigma_{pq}^\star} + \frac{\partial R_1}{\partial \tilde \Sigma_{pq}^\star}M_{mk}^{\star \star} \quad ; \quad \displaystyle \frac{\partial \tilde M_{mk}}{\partial M_{pqr}^\star} = R_1\displaystyle \frac{D \tilde M_{mk}^{\star \star}}{D \tilde M_{pqr}^\star} + \frac{\partial R_1}{\partial  \tilde M_{pqr}^\star}M_{mk}^{\star \star}.\nonumber
\end{equation}
Similarily, the fact that
\begin{equation}\label{eqn:ParamCrit2}
           \tilde M'_{ijk} = \displaystyle R_1\tilde {M_{ijk}^{\star \star}}' \quad \mbox {where} \quad R_2 = \displaystyle \sqrt{\frac{M_{II}}{M_{II}^{ \star \star}}} = \frac{\tilde \Sigma_{eq}}{\tilde \Sigma_{eq} + \displaystyle \frac{A_{II}}{5} (\tilde \Sigma_{eq}^\star-\tilde \Sigma_{eq})} = \frac{\tilde \Sigma_{eq}}{\mathcal {D}_2} \nonumber
\end{equation}
yields 
\begin{equation}
\label{eqn:ParamCrit2}
    \left\{
       \begin{array}{lll}
           \displaystyle \frac{\partial R_2}{\partial \tilde \Sigma_{pq}^\star} = \displaystyle \frac{1}{\mathcal D_2}\frac{D \tilde \Sigma_{eq}}{D \tilde \Sigma_{pq}^\star} - \frac{\tilde \Sigma_{eq}}{\mathcal D_2^2}\left( \frac{\partial \mathcal {D}_2}{\partial \tilde \Sigma_{pq}^\star} + \frac{\partial \mathcal D_2}{\partial \tilde \Sigma_{eq}}\frac{D \tilde \Sigma_{eq}}{D \tilde \Sigma_{pq}^\star}\right)\\
           {}&{}&{}\\
           \displaystyle \frac{\partial R_2}{\partial \tilde M_{pqr}^\star} = \displaystyle \frac{1}{\mathcal D_2}\frac{D \tilde \Sigma_{eq}}{D \tilde M_{pqr}^\star} - \frac{\tilde \Sigma_{eq}}{\mathcal {D}_2^2} \frac{\partial \mathcal {D}_2}{\partial \tilde \Sigma_{eq}}\frac{D \tilde \Sigma_{eq}}{D \tilde M_{pqr}^\star} \nonumber
       \end{array}
    \right.
\end{equation}
which is equivalent to
\begin{equation}\label{eqn:ParamCrit2}
 \displaystyle \frac{\partial \tilde M'_{ijk}}{\partial \Sigma_{pq}^\star} = R_2\displaystyle \frac{D \tilde {M_{ijk}^{\star \star}}'}{D \tilde \Sigma_{pq}^\star} + \frac{\partial R_2}{\partial \tilde \Sigma_{pq}^\star} {M_{ijk}^{\star \star}}' \quad ; \quad \displaystyle \frac{\partial \tilde M'_{ijk}}{\partial M_{pqr}^\star} = R_2\displaystyle \frac{D \tilde {M_{ijk}^{\star \star}}'}{D \tilde M_{pqr}^\star} + \frac{\partial R_2}{\partial \tilde M_{pqr}^\star} {M_{ijk}^{\star \star}}'.\nonumber
\end{equation}
We easily derive from these equations, the expressions of the derivatives: $ \displaystyle \frac{\partial \tilde M_{ijk}}{\partial \tilde \Sigma_{pq}^\star}$ and $ \displaystyle \frac{\partial \tilde M_{ijk}}{\partial \tilde M_{pqr}^\star}$.
We specify that
\begin{equation}
\label{eqn:DerivG5}
    \left\{
       \begin{array}{lll}
       \displaystyle \tilde M_{I} = \displaystyle \frac{\tilde M_{I}^{\star \star} \tilde \Sigma_{eq}^2}{ \left(\tilde \Sigma_{eq} + \displaystyle \frac{3\lambda+2\mu}{45 \mu}A_{I}(\tilde \Sigma_{eq}^\star-\tilde \Sigma_{eq})\right)^2}  = R_1^2 \tilde M_{I}^{\star \star}\\
           {}&{}&{}\\
        \tilde M_{II} = R_2^2 \tilde M_{II}^{\star \star}, \nonumber
       \end{array}
    \right.
\end{equation}
\begin{equation}
\label{eqn:DerivG5}
    \left\{
       \begin{array}{lll}
          \displaystyle \frac{\partial \tilde M_{I} }{\partial \tilde \Sigma_{pq}^\star} = \displaystyle 2R_1 \frac{\partial R_1}{\partial \tilde \Sigma_{pq}^\star} \tilde M_{I}^{\star \star} + R_1^2 \frac{D\tilde M_{I}^{\star \star}}{D\tilde \Sigma_{pq}^\star}\\
           {}&{}&{}\\
          \displaystyle \frac{\partial \tilde M_{II}}{\partial \tilde \Sigma_{pq}^\star} = \displaystyle 2R_2 \frac{\partial R_2}{\partial \tilde \Sigma_{pq}^\star} \tilde M_{II}^{\star \star} + R_2^2 \frac{D\tilde M_{II}^{\star \star}}{D\tilde \Sigma_{pq}^\star}\\
       \end{array}
    \right.
    \quad \mbox {with} \quad
    \left\{
       \begin{array}{lll}
         \displaystyle \frac{D \tilde M_{I}^{\star \star}}{D \tilde \Sigma_{pq}^\star} = \displaystyle \frac{\partial \tilde M_{I}^{\star \star}}{\partial \tilde \Sigma_{pq}^\star} + \frac{\partial \tilde M_{I}^{\star \star}}{\partial \tilde \Sigma_{eq}}\frac{D \tilde \Sigma_{eq}}{D \tilde \Sigma_{pq}^\star}\\
           {}&{}&{}\\
         \displaystyle \frac{D \tilde {M_{II}^{\star \star}}'}{D \tilde \Sigma_{pq}^\star} = \displaystyle \frac{\partial \tilde {M_{II}^{\star \star}}'}{\partial \tilde \Sigma_{pq}^\star} + \frac{\partial \tilde {M_{II}^{\star \star}}'}{\partial \tilde \Sigma_{eq}}\frac{D \tilde \Sigma_{eq}}{D \tilde \Sigma_{pq}^\star}.\nonumber
       \end{array}
    \right.
\end{equation}
%%
%%%%%%%%%%%%%%%%%%%%%%%%%
\subsubsection{Introduction of work hardening and transition to non-reduced stresses}
\label{deriv27}
The equation giving $ \Delta \bar \varepsilon$, written with ``reduced'' constraints and moments, is:
\begin{eqnarray}
\label{eqn:DerivG10}
    (1 - f) \frac{\Delta \bar \varepsilon}{\bar \Sigma} &=& \tilde \Sigma_{eq}\frac{\tilde \Sigma_{eq}^\star - \tilde \Sigma_{eq}}{3 \mu} + 3 \tilde \Sigma_m \frac{\tilde \Sigma_m^\star - \tilde \Sigma_m}{3 \lambda + 2 \mu} + \nonumber\\
 &+& \tilde M_{ijk}' \frac{\tilde {M_{ijk}^\star}' - \tilde M'_{ijk}}{2 \mu b^2/5} + 3 \tilde M_{mk} \frac{\tilde M_{mk}^\star - \tilde M_{mk}}{(3 \lambda + 2 \mu)b^2/5} \nonumber
\end{eqnarray}
which can also be written in the form:
\begin{eqnarray}
\label{eqn:DerivG11}
 (1 - f) \frac{\Delta \bar \varepsilon}{\bar \Sigma} &=& \tilde \Sigma'_{ij}\frac{\tilde {\Sigma_{ij}^\star}' - \tilde \Sigma'_{ij}}{2 \mu} + 3 \tilde \Sigma_m \frac{\tilde \Sigma_m^\star - \tilde \Sigma_m}{3 \lambda + 2 \mu} + \nonumber\\
& + & \tilde M'_{ijk} \frac{\tilde {M_{ijk}^\star}' - \tilde M'_{ijk}}{2 \mu b^2/5} + 3 \tilde M_{mk} \frac{\tilde M_{mk}^\star - \tilde M_{mk}}{(3 \lambda + 2 \mu)b^2/5} \nonumber
\end{eqnarray}
Differentiating the last expression by considering small variations of $ \tilde \Sigma_{pq}^\star$ and $ \tilde M_{pqr}^\star$, we get
\begin{eqnarray}
\label{eqn:DerivG11}
    D \left[ (1 - f) \frac {\Delta \bar \varepsilon}{\bar \Sigma} \right] &=&  (1 - f) \left[ \frac{D \Delta \bar \varepsilon}{\bar \Sigma} - \frac{\Delta \bar \varepsilon}{\bar \Sigma^2} D \bar \Sigma \right] \nonumber\\
&=& \frac{1 - f}{\bar \Sigma} ( \frac{1}{h} - \frac {\Delta \bar \varepsilon}{\bar \Sigma})D \bar \Sigma \nonumber
\end{eqnarray}
where
\begin{eqnarray}
\label{eqn:DerivG11}
   h = \frac{D \bar \Sigma}{D \bar \varepsilon}. \nonumber
\end{eqnarray}
Thus, we obtain
\begin{eqnarray}
\label{eqn:DerivG12}
    \frac{1 - f}{\bar \Sigma}( \frac{1}{h} - \frac{\Delta \bar \varepsilon}{\bar \Sigma}) D \bar \Sigma &=& D \tilde \Sigma'_{ij}\frac{\tilde {\Sigma_{ij}^\star}' - \tilde \Sigma'_{ij}}{2 \mu} + 3 D \tilde \Sigma_m \frac{\tilde \Sigma_m^\star - \tilde \Sigma_m}{3 \lambda + 2 \mu} \nonumber\\
     &+& D \tilde M'_{ijk} \frac{\tilde {M_{ijk}^\star}' - \tilde M'_{ijk}}{2 \mu b^2/5} + 3 D \tilde M_{mk} \frac{\tilde M_{mk}^\star - \tilde M_{mk}}{(3 \lambda + 2 \mu)b^2/5} \nonumber\\
     & + & \tilde \Sigma_{ij}' \frac{ D \tilde {\Sigma_{ij}^\star}' - D \tilde \Sigma_{ij}'}{2 \mu} +  \frac { 3 \tilde \Sigma_m }{ 3 \lambda + 2 \mu } (  D \tilde \Sigma_{m}^\star - D \tilde \Sigma_m ) \nonumber\\
     &+& \tilde M'_{ijk} \frac{D \tilde {M_{ijk}^\star}' - D \tilde M'_{ijk}}{2 \mu b^2/5} + 3 \tilde M_{mk} \frac{D \tilde M_{mk}^\star - D \tilde M_{mk}}{(3 \lambda + 2 \mu)b^2/5} \nonumber
\end{eqnarray}
The sum of the first four terms is zero.
Indeed, this sum can be written as
\begin{eqnarray}
\label{eqn:DerivG11}
    D \tilde \Sigma'_{ij}  \Delta \tilde {\varepsilon_{ij}^p}' + 3\mbox {d} \tilde \Sigma_m  \Delta \tilde \varepsilon_m^p + D \tilde M'_{ijk} \Delta {(\nabla \tilde w)_{ijk}^p}' + 3 D \tilde M_{mk} (\nabla \tilde w)_{mk}^p = D \tilde \Sigma_{ij}  \Delta \tilde \varepsilon_{ij}^p + D \tilde M_{ijk} \Delta (\nabla \tilde w)_{ijk}^p.\nonumber
\end{eqnarray}
Yet, the plasticity criterion of Gologanu {\it et al.} can be written in the form
\begin{eqnarray}
\label{eqn:DerivG12}
    \tilde \Phi (\tilde {\bf \Sigma}, \tilde {\bf M}) = 0 \quad \Rightarrow \quad \frac{\partial \tilde \Phi}{\partial \tilde \Sigma_{ij}} D \tilde \Sigma_{ij} + \frac{\partial \tilde \Phi}{\partial \tilde M_{ijk}} D \tilde M_{ijk} = 0 \quad \mbox {and} \quad \Delta  \tilde \epsilon_{ij}^p = \Delta \tilde \eta \frac{\partial \tilde \Phi}{\partial \tilde \Sigma_{ij}} \quad , \quad \nonumber\\
   \Delta (\nabla \tilde w)_{ijk}^p = \Delta \tilde \eta \left( \frac{\partial \tilde \Phi}{\partial \tilde M_{ijk}} + \delta_{ik}\tilde U_j + \delta_{jk}\tilde U_i \right). \nonumber
\end{eqnarray}
As a consequence,
\begin{eqnarray}
\label{eqn:DerivG12}
   \Delta \tilde \varepsilon_{ij}^p D \tilde \Sigma_{ij}   +  \Delta (\nabla \tilde w)_{ijk}^p D \tilde M_{ijk} =  \Delta \tilde \eta \left( \frac{\partial \tilde \Phi}{\partial \tilde \Sigma_{ij}} D \tilde \Sigma_{ij} + \frac{\partial \tilde \Phi}{\partial \tilde M_{ijk}} D \tilde M_{ijk} +  \delta_{ik}\tilde U_j D \tilde M_{ijk} + \delta_{jk}\tilde U_i D \tilde M_{ijk}\right) = 0 \nonumber
\end{eqnarray}
because
\begin{eqnarray}
\label{eqn:DerivG12}
D \tilde M_{ijj} = 0.\nonumber
\end{eqnarray}
\subsubsection{The derivatives of $ \bar \Sigma$ par rapport aux $ \tilde \Sigma_{pq}^\star$ et aux $ \tilde M_{pqr}^\star$}.
\label{deriv}
We get
\begin{equation}
\label{eqn:DerivG12}
    \left\{
       \begin{array}{lll}
    \displaystyle \frac{1 - f}{\bar \Sigma} \left( \frac{1}{h} - \frac{\Delta \bar \varepsilon}{\bar \Sigma} \right) \frac{\partial \bar \Sigma}{\partial \tilde \Sigma_{pq}^\star} &=&  \displaystyle \frac{\tilde \Sigma'_{ij}}{2 \mu} \left(\frac{\partial \tilde {\Sigma_{ij}^\star}'}{\partial \tilde \Sigma_{pq}^\star} - \frac{\partial \tilde \Sigma'_{ij}}{\partial \tilde \Sigma_{pq}^\star} \right)  + \frac{3 \tilde \Sigma_m}{3 \lambda + 2 \mu} \left( \frac{\partial \tilde \Sigma_m^\star}{\partial \tilde {\Sigma_{pq}^\star}} - \frac{\partial \tilde \Sigma_m}{\partial \tilde \Sigma_{pq}^\star} \right) \nonumber\\
     &-&  \displaystyle \frac {\tilde M'_{ijk}}{2 \mu b^2/5} \frac{\partial \tilde {M_{ijk}}'}{\tilde \partial \Sigma_{pq}^\star} - 3 \frac {\tilde M_{mk}} {(3 \lambda + 2 \mu)b^2/5} \frac{\partial \tilde M_{mk}}{\tilde \partial \Sigma_{pq}^\star} \nonumber\\
    {}&{}&{}\nonumber\\
     &=&  \displaystyle \frac{\tilde \Sigma'_{ij}}{2 \mu} \left[ \frac{1}{2}(\delta_{ip}\delta_{jq} + \delta_{iq}\delta_{jp}) - \frac{1}{3}\delta_{ij}\delta_{pq}\right] - \frac{1}{6 \mu} \frac{\mbox D \tilde \Sigma_{eq}^2}{\mbox D \tilde \Sigma_{pq}^\star} \nonumber\\
     &+&  \displaystyle \frac{3 \tilde \Sigma_m}{3 \lambda + 2 \mu} \frac{1}{3} \delta_{pq} - 3 \frac{\tilde \Sigma_m}{3\lambda + 2 \mu} \frac{ \partial \tilde \Sigma_m}{\partial \tilde \Sigma_{pq}^\star}  \nonumber\\
     &-&  \displaystyle \frac{1}{6 \mu b^2/5} \frac{\partial \tilde M_{II}}{\partial \tilde \Sigma_{pq}^\star} - 3 \frac {\tilde M_{mk}} {(3 \lambda + 2 \mu)b^2/5} \frac{\partial \tilde M_{mk}}{\partial \tilde \Sigma_{pq}^\star} \nonumber\\
    {}&{}&{}\nonumber\\
     &=&  \displaystyle \frac{\tilde \Sigma'_{pq}}{2 \mu} - \frac{\tilde \Sigma_{eq}}{3 \mu}\frac{\mbox D \tilde \Sigma_{eq}}{\mbox D \tilde \Sigma_{pq}^\star} + \frac{\tilde \Sigma_m}{3 \lambda + 2 \mu} \delta_{pq} - 3 \frac{\tilde \Sigma_m}{3\lambda + 2 \mu} \frac{ \partial \tilde \Sigma_m}{\partial \tilde \Sigma_{pq}^\star}  \nonumber\\
     &-&  \displaystyle \frac{1}{6 \mu b^2/5} \frac{\partial \tilde M_{II}}{\partial \tilde \Sigma_{pq}^\star} - \frac {3}{2(3 \lambda + 2 \mu)b^2/5}\frac{\partial \tilde M_{I}}{\partial \tilde \Sigma_{pq}^\star}. \nonumber
       \end{array}
    \right.
\end{equation} 
Similarily,
\begin{equation}
\label{eqn:DerivG12}
    \left\{
       \begin{array}{lll}
    \displaystyle \frac{1 - f}{\bar \Sigma} \left( \frac{1}{h} - \frac{\Delta \bar \varepsilon}{\bar \Sigma} \right) \frac{\partial \bar \Sigma}{\partial \tilde M_{pqr}^\star} &=&  - \displaystyle \frac{\tilde \Sigma'_{ij}}{2 \mu} \frac{\partial \tilde \Sigma'_{ij}}{\partial \tilde M_{pqr}^\star}  - \frac{3 \tilde \Sigma_m}{3 \lambda + 2 \mu}  \frac{\partial \tilde \Sigma_m}{\partial \tilde M_{pqr}^\star} \nonumber \\
 & + & \displaystyle \frac{\tilde M'_{ijk}}{2 \mu b^2 / 5} \left[ \frac{1}{2}(\delta_{ip}\delta_{jq} + \delta_{iq}\delta_{jp}) - \frac{1}{3}\delta_{ij}\delta_{pq}\right]\delta_{kr} \nonumber\\
     &-&  \displaystyle \frac {\tilde M'_{ijk}}{2 \mu b^2/5} \frac{\partial \tilde {M_{ijk}}'}{\partial \tilde M_{pqr}^\star} + 3 \frac {\tilde M_{mk}} {(3 \lambda + 2 \mu) b^2/5}\frac{1}{3} \delta_{pq}\delta_{kr} \nonumber\\
& - & 3  \displaystyle \frac {\tilde M_{mk}} {(3 \lambda + 2 \mu)b^2/5} \frac{\partial \tilde M_{mk}}{\partial \tilde M_{pqr}^\star} \nonumber\\
    {}&{}&{}\nonumber\\
     &=& - \displaystyle \frac{\tilde \Sigma_{eq}}{3 \mu}\frac{\mbox D \tilde \Sigma_{eq}}{D \tilde M_{pqr}^\star} + \frac{\tilde M_{mr}}{(3 \lambda + 2 \mu)b^2/5} \delta_{pq} - 3 \frac{\tilde \Sigma_m}{3\lambda + 2 \mu} \frac{ \partial \tilde \Sigma_m}{\partial \tilde M_{pqr}^\star} \nonumber\\
     &-& \displaystyle \frac{1}{6 \mu b^2/5} \frac{\partial \tilde M_{II}}{\partial \tilde M_{pqr}^\star} - \frac {3}{2 (3 \lambda + 2 \mu)b^2/5}\frac{\partial \tilde M_{I}}{\partial \tilde M_{pqr}^\star} + \frac{\tilde M'_{pqr}}{2 \mu b^2/5}. \nonumber
       \end{array}
    \right.
\end{equation}
%%
%%%%%%%%%%%%%%%%%%%%%%%%%
\subsubsection{The derivatives of $ \bar \Sigma$ with respect $ \Sigma_{pq}^\star$ and $ M_{pqr}^\star$}
\label{deriv}
%%%%%%%%%%%%%%%%%%%%%%%%%
%%
Finally, we determine the expression of the non-normalized stresses and moments:
\begin{eqnarray}\label{eqn:DerivG12}
    \bar \Sigma  &\equiv& \bar \Sigma (\tilde {\bf \Sigma}^\star, {\bf M}^\star) \equiv \bar \Sigma (\frac{\tilde {\bf \Sigma}^\star}{\bar \Sigma}, \frac{{\bf M}^\star}{\bar \Sigma}) \quad \quad \quad \quad \quad \quad \quad \quad \quad \quad \quad \quad \quad \quad \quad \quad \nonumber \\ 
    &\Rightarrow& \quad D \bar \Sigma = \frac{\partial \bar \Sigma}{\partial \tilde \Sigma_{pq}^\star} \left( \frac{D \Sigma_{pq}^\star}{\bar \Sigma} - \frac{\Sigma_{pq}^\star}{\bar \Sigma^2}{D \bar \Sigma} \right) + \frac{\partial \bar \Sigma}{\partial \tilde M_{pqr}^\star} \left( \frac{D M_{pqr}^\star}{\bar \Sigma} - \frac{M_{pqr}^\star}{\bar \Sigma^2}{D \bar \Sigma} \right)\quad \quad \quad \nonumber\\
    &\Rightarrow& \left( 1 + \frac{1}{\bar \Sigma^2}\frac{\partial \bar \Sigma}{\partial \tilde \Sigma_{pq}^\star} \Sigma_{pq}^\star + \frac{1}{\bar \Sigma^2}\frac{\partial \bar \Sigma}{\partial \tilde M_{pqr}^\star} M_{pqr}^\star \right) D \bar \Sigma = \frac{1}{\bar \Sigma} \frac{\partial \bar \Sigma}{\partial \tilde \Sigma_{pq}^\star} D \Sigma_{pq}^\star + \frac{1}{\bar \Sigma} \frac{\partial \bar \Sigma}{\partial \tilde M_{pqr}^\star} D M_{pqr}^\star. \nonumber
\end{eqnarray}
Finally, we obtain
\begin{equation}
\label{eqn:DerivG5}
    \left\{
       \begin{array}{lll}
       \displaystyle \frac{\partial \bar \Sigma}{\partial \Sigma_{pq}^\star} = \displaystyle \frac{\displaystyle \frac{1}{\bar \Sigma} \displaystyle\frac{\partial \bar \Sigma}{\partial \tilde \Sigma_{pq}^\star}}{1 + \displaystyle\frac{1}{\bar \Sigma^2}\frac{\partial \bar \Sigma}{\partial \Sigma_{ij}^\star} \Sigma_{ij}^\star + \displaystyle\frac{1}{\bar \Sigma^2}\frac{\partial \bar \Sigma}{\partial \partial M_{ijk}^\star} M_{ijk}^\star}\\
           {}&{}&{}\\
       \displaystyle \frac{\partial \bar \Sigma}{\partial M_{pqr}^\star} = \displaystyle \frac{\displaystyle \frac{1}{\bar \Sigma} \displaystyle\frac{\partial \bar \Sigma}{\partial \tilde M_{pqr}^\star}}{1 + \displaystyle\frac{1}{\bar \Sigma^2}\displaystyle\frac{\partial \bar \Sigma}{\partial \tilde \Sigma_{ij}^\star} \Sigma_{ij}^\star + \displaystyle\frac{1}{\bar \Sigma^2}\frac{\partial \bar \Sigma}{\partial \tilde M_{ijk}^\star} M_{ijk}^\star}. \nonumber
       \end{array}
    \right.
\end{equation}
We had $\Sigma_{ij} \equiv \bar \Sigma (\bs^\star , {\bf M}^\star) \tilde \Sigma_{ij} (\displaystyle \frac{\bs^\star}{ \bar \Sigma} , \displaystyle \frac{{\bf M}^\star}{\bar \Sigma})$; and it follows that
\begin{eqnarray}
\label{eqn:DerivG12}
 D \Sigma_{ij} &=& D \bar \Sigma \tilde \Sigma_{ij} + \bar \Sigma \frac{\partial \tilde \Sigma_{ij}}{\partial \tilde \Sigma_{pq}^\star} \left( \frac{ D \Sigma_{pq}^\star}{\bar \Sigma} - \frac { \Sigma_{pq}^\star}{{\bar \Sigma}^2} D \bar \Sigma \right) + \bar \Sigma \frac{\partial \tilde \Sigma_{ij}}{\partial \tilde M_{pqr}^\star} \left( \frac{D M_{pqr}^\star}{\bar \Sigma} - \frac{M_{pqr}^\star}{{\bar \Sigma}^2} D \bar \Sigma \right) \nonumber\\
     &=& \left( \tilde \Sigma_{ij} - \frac{\partial \tilde \Sigma_{ij}}{\partial \tilde \Sigma_{hk}^\star}\frac{\Sigma_{hk}^\star}{\bar \Sigma} - \frac{\partial \tilde \Sigma_{ij}}{\partial \tilde M_{hkl}^\star}\frac{M_{hkl}^\star}{\bar \Sigma} \right) \left( \frac{\partial \bar \Sigma}{\partial \Sigma_{pq}^\star} D \Sigma_{pq}^\star + \frac{\partial \bar \Sigma}{\partial \tilde M_{pqr}^\star} D M_{pqr}^\star\right) \nonumber\\
     &+& \frac{\partial \tilde \Sigma_{ij}}{\partial \tilde \Sigma_{pq}^\star} D \Sigma_{pq}^\star +  \frac{\partial \tilde \Sigma_{ij}}{\partial \tilde M_{pqr}^\star} D M_{pqr}^\star. \nonumber
\end{eqnarray}
Thereafter,
\begin{equation}\label{eqn:DerivG5}
    \left\{
       \begin{array}{lll}
       \displaystyle \frac{\partial \Sigma_{ij}}{\partial \Sigma_{pq}^\star} = \left( \tilde \Sigma_{ij} - \frac{\partial \tilde \Sigma_{ij}}{\partial \tilde \Sigma_{hk}^\star}\tilde \Sigma_{hk}^\star - \frac{\partial \tilde \Sigma_{ij}}{\partial \tilde M_{hkl}^\star}\tilde M_{hkl}^\star \right)\frac{\partial \bar \Sigma}{\partial \Sigma_{pq}^\star} + \frac{\partial \tilde \Sigma_{ij}}{\partial \tilde \Sigma_{pq}^\star}\\
           {}&{}&{}\\
       \displaystyle \frac{\partial \Sigma_{ij}}{\partial M_{pqr}^\star} = \left( \tilde \Sigma_{ij} - \frac{\partial \tilde \Sigma_{ij}}{\partial \tilde \Sigma_{hk}^\star}\tilde \Sigma_{hk}^\star - \frac{\partial \tilde \Sigma_{ij}}{\partial \tilde M_{hkl}^\star} \tilde M_{hkl}^\star \right)\frac{\partial \bar \Sigma}{\partial M_{pqr}^\star} + \frac{\partial \tilde \Sigma_{ij}}{\partial \tilde M_{pqr}^\star}. \nonumber
       \end{array}
    \right.
\end{equation}
Similarily, $M_{ijk} \equiv \bar \Sigma (\bs^\star , {\bf M}^\star) \tilde M_{ijk} (\displaystyle \frac{\bs ^\star}{\bar \Sigma} , \displaystyle \frac{{\bf M}^\star}{\bar \Sigma})$; we deduce that
\begin{equation}
\label{eqn:DerivG5}
    \left\{
       \begin{array}{lll}
       \displaystyle \frac{\partial M_{ijk}}{\partial \Sigma_{pq}^\star} = \left( \tilde M_{ijk} - \frac{\partial \tilde M_{ijk}}{\partial \tilde \Sigma_{lm}^\star}\tilde \Sigma_{lm}^\star - \frac{\partial \tilde M_{ijk}}{\partial \tilde M_{lmn}^\star}\tilde M_{lmn}^\star \right)\frac{\partial \bar \Sigma}{\partial \Sigma_{pq}^\star} + \frac{\partial \tilde M_{ijk}}{\partial \tilde \Sigma_{pq}^\star}\\
           {}&{}&{}\\
       \displaystyle \frac{\partial M_{ijk}}{\partial M_{pqr}^\star} = \left( \tilde M_{ijk} - \frac{\partial \tilde M_{ijk}}{\partial \tilde \Sigma_{lm}^\star}\tilde \Sigma_{lm}^\star - \frac{\partial \tilde M_{ijk}}{\partial \tilde M_{lmn}^\star} \tilde M_{lmn}^\star \right)\frac{\partial \bar \Sigma}{\partial M_{pqr}^\star} + \frac{\partial \tilde M_{ijk}}{\partial \tilde M_{pqr}^\star}. \nonumber
       \end{array}
    \right.
\end{equation}
%%
%%%%%%%%%%%%%%%%%%%%%%%%%
\subsubsection{The derivatives of $ \Sigma_{ij}$ and $ M_{ijk}$ with respect to $ \Delta \varepsilon_{pq}$ and $ \Delta (\nabla {W})_{pqr}$}
\label{deriv}
%%%%%%%%%%%%%%%%%%%%%%%%%
We obtain
\begin{equation}\label{eqn:DerivG5}
    \left\{
       \begin{array}{lll}
       \displaystyle \frac{\partial \Sigma_{ij}}{\partial \Delta \varepsilon_{pq}} = \displaystyle \frac{\partial \Sigma_{ij}}{\partial \Sigma_{lm}^\star}\frac{\partial \Sigma_{lm}^\star}{\partial \Delta \varepsilon_{pq}} \quad ; \quad \frac{\partial \Sigma_{ij}}{\Delta (\nabla {W})_{pqr}} = \frac{\partial \Sigma_{ij}}{\partial M_{lmn}^\star} \frac{\partial M_{lmn}^\star}{\partial \Delta (\nabla {W})_{pqr}}\\
        {}&{}&{}\\
       \displaystyle \frac{\partial M_{ijk}}{\partial \Delta \varepsilon_{pq}} = \displaystyle \frac{\partial M_{ijk}}{\partial \Sigma_{lm}^\star}\frac{\partial \Sigma_{lm}^\star}{\partial \Delta \varepsilon_{pq}} \quad ; \quad \frac{\partial M_{ijk}}{\Delta (\nabla {W})_{pqr}} = \frac{\partial M_{ijk}}{\partial M_{lmn}^\star} \frac{\partial M_{lmn}^\star}{\Delta (\nabla {W})_{pqr}}.\nonumber
       \end{array}
    \right.
\end{equation}
\section{Assessment of the approach}
\label{sec:StiffMatrix}
\subsection{Simulations of axisymmetric notched specimen}
\label{sec:StiffMatrix}
We have implemented the developed tangent stiffness matrix for the GLPD model using Systus$^{\circledR}$ finite element software developed by ESI group. To test the procedure we simulate 
a tensile test on an axisymmetric notched specimen in A508 Cl.3 steel, for which the mechanical fields (stresses and deformations) are homogeneous in the ligament of the specimen.
$\\$

In each of the simulations, we use the value of the Tvergaard parameter $ q = 4/e = 1.47$ determined by Perrin and Leblond \textcolor{blue} { \cite{PL90} } by a ``differential' scheme '.
The  A508 Cl.3 steel is used in the design of reactor vessels for nuclear power plants.
For this steel, the Young's modulus is worth $ {E} = 203 000$ MPa, the Poisson's ratio $\nu = 0.3$, and the initial yield stress in simple tension $ \sigma_y = 450 $ MPa.
The value of the initial porosity $ f_0$ of this steel comes from the study of its chemical composition and its inclusionary state \textcolor{blue} { \cite{RM81} }. The value of $ f_0$ is deduced from the sulfur and manganese content of the material and the average dimensions of the inclusions. 
The value found is $0.00016$.
To this damage parameter are added two others, $f_{c} $ and $ \delta $, respectively representing the critical porosity at the start of coalescence and the accelerating factor of the growth of the cavities whose values can be adjusted from one simulation to another.
$\\$

\begin{figure}[hbt!]
\begin{center}
\includegraphics[width=5.5 cm, height= 12 cm, angle=0]{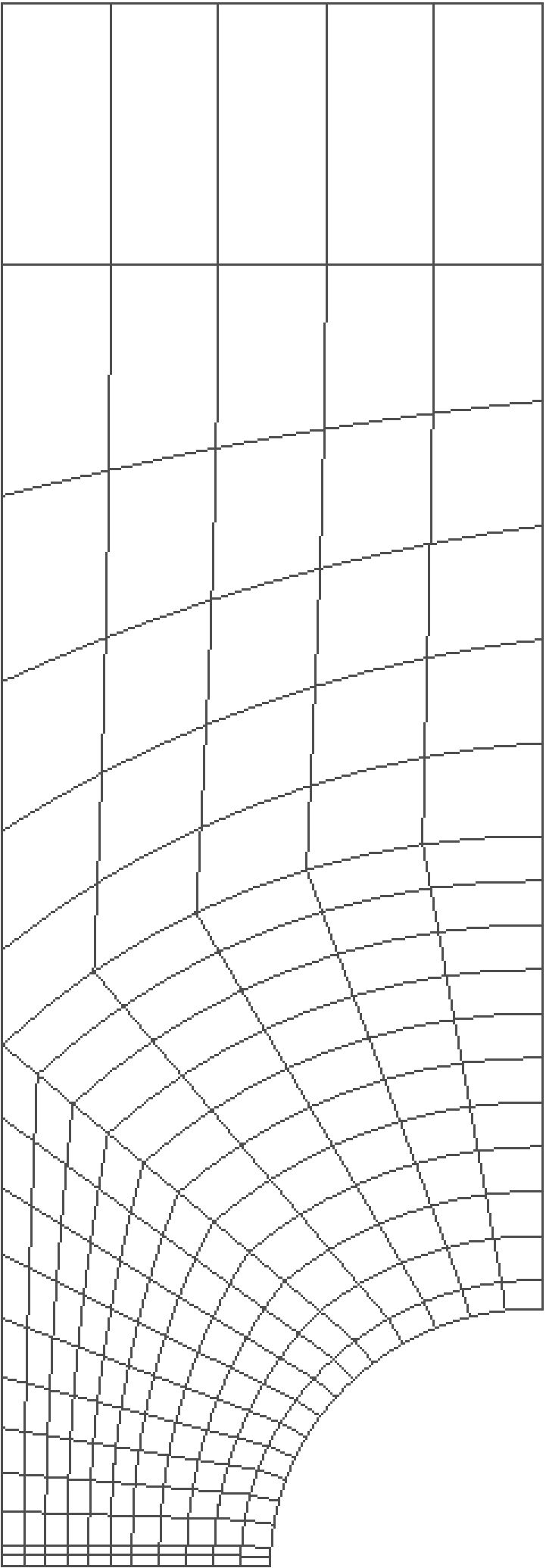}
   \caption{\textit{Mesh of the AE05 Specimen}}
\label{fig:002}
\end{center}
\end{figure}

The notch radius of this specimen is $5 mm$.
The specimen used is cylindrical.
The symmetry and axisymmetry of the problem are exploited to model only a quarter of a longitudinal section of the specimen.
The nodes located in the minimum section are axially locked in the direction ${ Y}$, the loading of the specimen is a displacement applied to all the nodes of the upper surface of the discretized sample portion, varying linearly over time.
The numerical calculations are carried out with meshes whose dimensions are worth $0.4 \times 0.2  mm^2$ in the region of the ligament, for a characteristic distance $ b$ of $10^{-5 }$ $\mu$m.
The mesh of this specimen is given in Figure \ref{fig:002}.
\begin{figure}[hbt!]
\begin{center}
  \includegraphics[width=8 cm, height = 8 cm, angle=0]{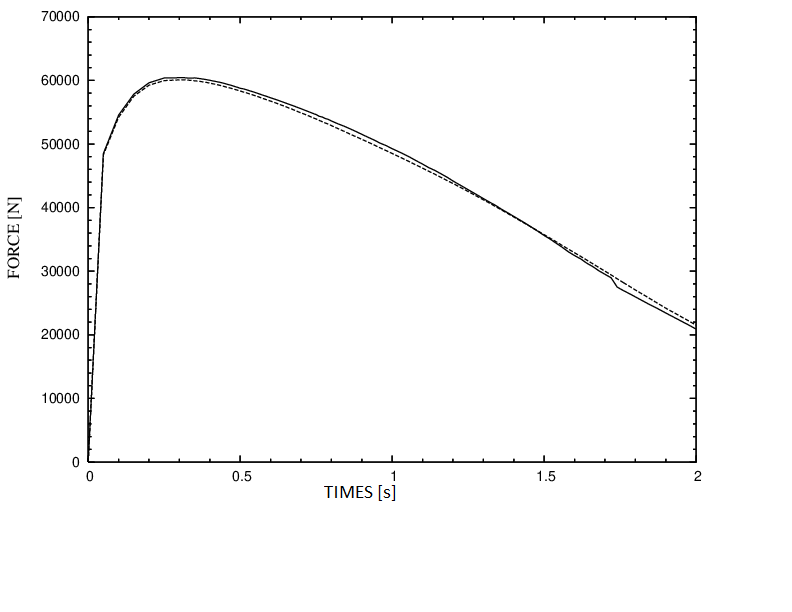}
   \caption{\textit{ Force-displacement curves full line time: Gurson model; solid line: GLPD model}}
\label{fig:003}
\end{center}
\end{figure}
$\\$

Figure \ref{fig:003} shows the force-displacement curves obtained with the Gurson model and the GLPD model.
In this specific case, we do not make a comparison with experiments because we do not have experimental results for this type of specimen.
On the other hand, we compare between them the curves obtained with the local Gurson models and the GLPD model for a characteristic distance $ b = 0.00001 \mu$m .
These siimulations indicate that for values of the characteristic distance tending toward $0$, we find the load curve obtained with the local Gurson's model.
These comparisons seem to indicate that the numerical results are going in the right direction and allow us to continue the tests on specimen geometries where the stress and strain gradients are large.
Note that quadratic convergence of elastoplastic iterations were obtained without using the stiffness matrix relations we discussed in the previous sections. The situation will be different with the simulations
of ductile fracture of pre-cracked specimens.   
\subsection{Simulations of pre-cracked specimens}
\label{sec:StiffMatrix}
\begin{figure}[hbt!]
\begin{center}
  \includegraphics[width=5.5cm, height = 12cm, angle=0]{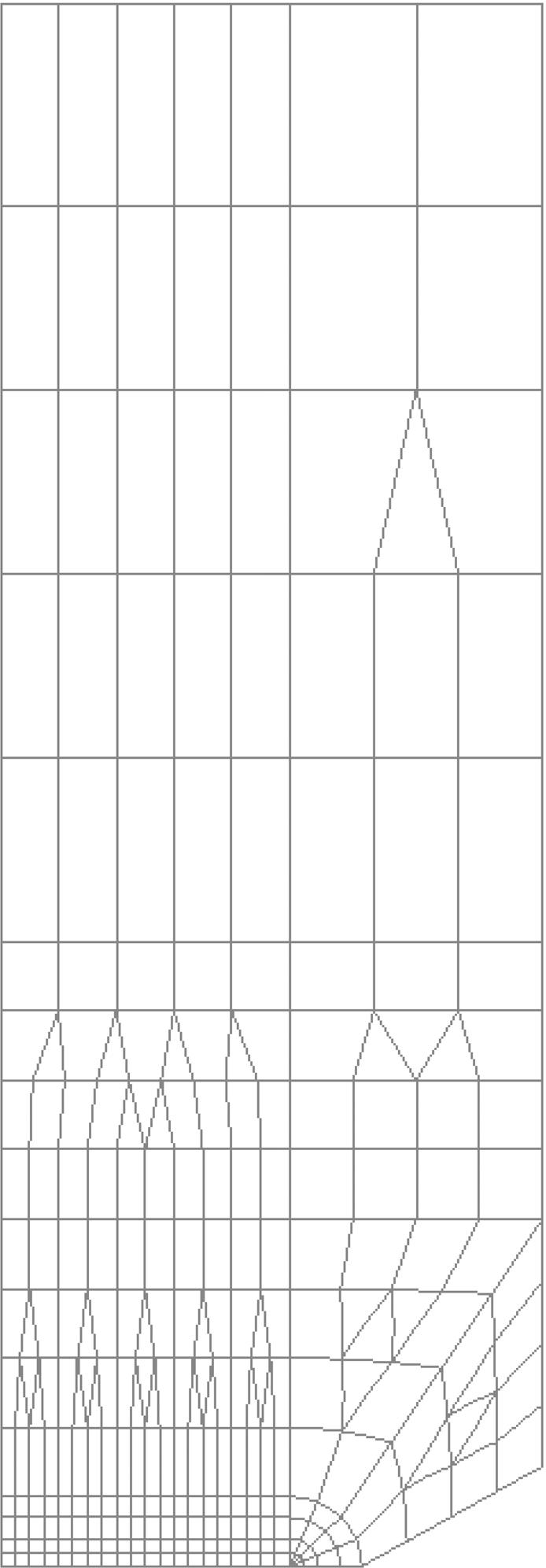}
   \caption{\textit{Mesh of the TA15 specimen}}
\label{fig:005}
\end{center}
\end{figure}

We now consider the axisymmetric pre-cracked specimen TA15 in A508 Cl.3 steel, for which we have experimental results.
We take advantage of the conditions of symmetry and axisymmetry by modeling only a quarter (see Figure \ref{fig:005} ) of a longitudinal section of the specimen with again quadrangle elements (8 nodes and 4 Gaussian points per sub-integrated element).
The few triangular elements, unavoidable in the automatic generation of the mesh, have no effect on the result because they are locatedoutside the sensitive areas of this problem (area of the crack tip and region ahead of the crack tip).
%The few triangular elements, unavoidable in the automatic generation of the mesh, have no effect on the result because they are located outside the ``sensitive'' areas of this problem (area of ​​the crack tip and region ahead of this tip).
%
The values of the half-height and the radius of the TA15 specimen are 22.5 $mm$ and 15 $ mm $. 
%%
%The values ​​of the half-height and the radius of the TA15 specimen are 22.5 $ mm $ and 15 $ mm $.
%
The half-angle of the opening is worth $15^{\circ}$; the depth of the central V-shaped notch is 3.88 $ mm $.
%The half-angle of opening is worth $15^{\circ}$; the depth of the central V-shaped notch is 3.88 $ mm $.
%
The pre-cracking radii is taken equal to 4 $mm$.
%The pre-cracking radii are taken equal to 4 $ mm$.
%%
$\\$

\begin{figure}[hbt!]
\begin{center}
  \includegraphics[width=4 cm, height = 4 cm, angle=0]{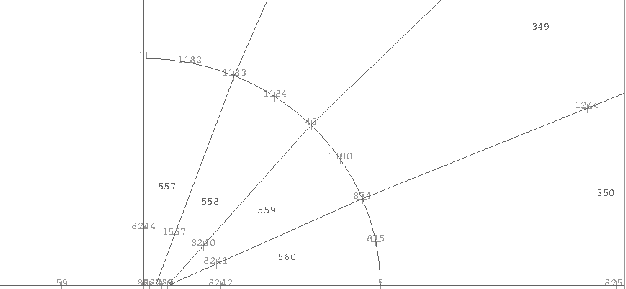}
   \caption{\textit{Mesh at the tip of the crack with 4 meshes whose nodes are pushed back to the quarter}}
\label{fig:004}
\end{center}
\end{figure}

Behind the crack, the discretization radiates with an angular sector divided into 4 to properly represent the significant stress and strain gradients in this area.
Among the four meshes adjoining the crack tip whose intermediate nodes are pushed back to the quarter, we distinguish two quasi-degenerate quadrilaterals and two triangles.
The role of the first meshes is to enable the representation of the blunting; the quasi-merged nodes deviate during the deformation.
In front of the crack, the discretization comprises identical square elements, see Figure \ref{fig:004}.
\begin{figure}[hbt!]
\begin{center}
  \includegraphics[width=6.5cm, height = 5.5 cm, angle=0]{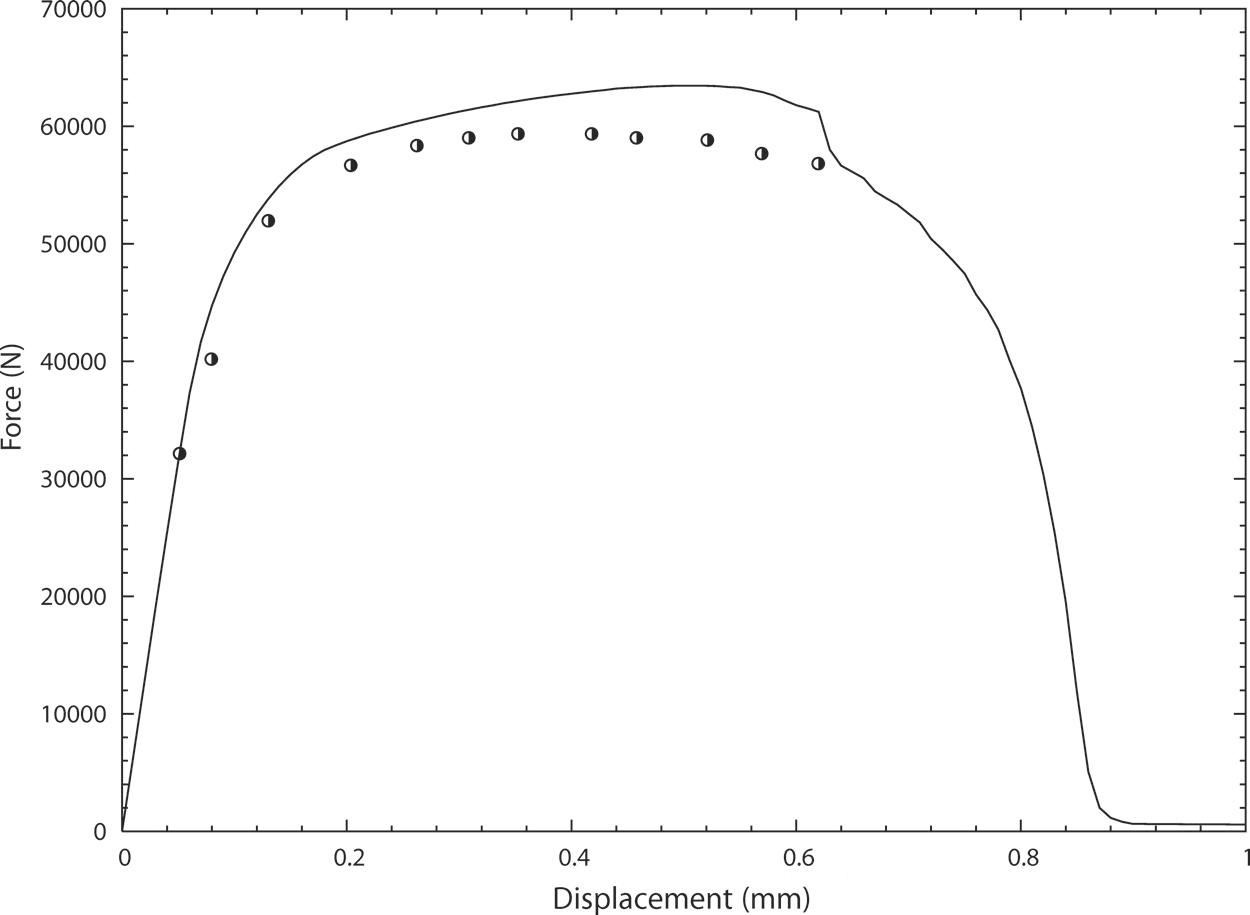}
   \caption{\textit{Stress strain curve foir the TA15 specimen with the original GLPD algorithm}}
\label{fig:006}
\end{center}
\end{figure}
\begin{figure}[hbt!]
\begin{center}
  \includegraphics[width=6.5cm, height = 5.5 cm, angle=0]{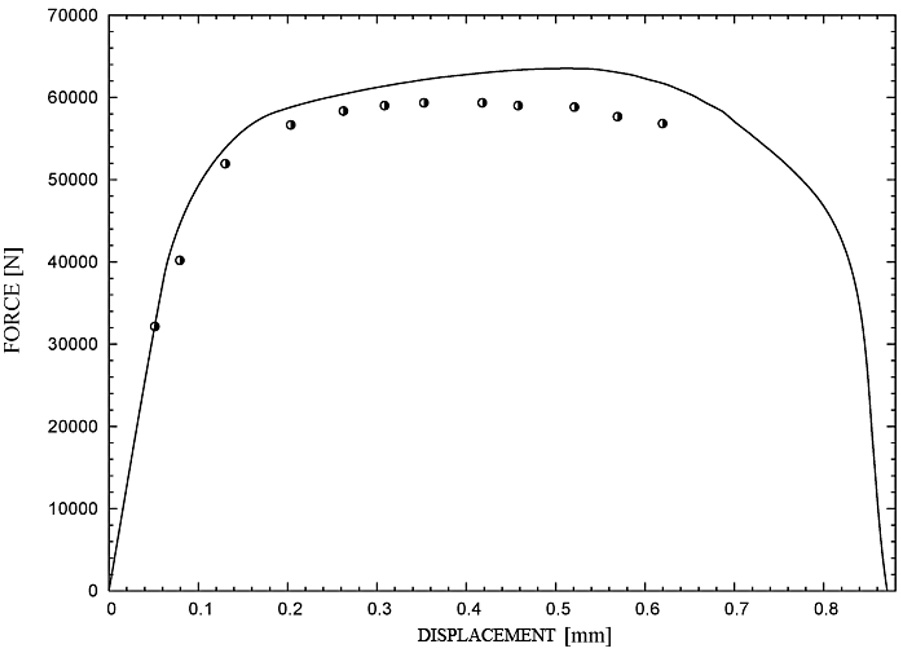}
   \caption{\textit{Stress strain curve foir the TA15 specimen with the modified GLPD algorithm}}
\label{fig:007}
\end{center}
\end{figure}
The numerical results are obtained using the GLPD model with a characteristic distance $ b=550\, \mu$m for a mesh size of $ 200\, \mu$m in the region of the ligament (for both specimens).
$\\$

The simulations on the TA15 specimen could not be completed.
Indeed, it was impossible to obtain the convergence of calculations beyond a certain level of loading as indicated in Figure\ref{fig:006} which displays the experimental and numerical force-displacement curves.
The numerical experiment shows that the use of the stifness tangent matrix allows a faster convergence (in number of iterations, otherwise in computing time) in cases where a
classical BFGS method also leads to convergence, but brings no improvement in the case where the BFGS method
leads to elasto-plastic convergence issues.
$\\$

An explicit variant of the elasto-plastic algorithm that does not present these convergence issues ( see Fgiure \ref{fig:007} ) consists in fixing the values of the increments of plastic deformation and of its gradient to those of the previous time step.
%%
%An explicit variant of the elasto-plastic algorithm that does not present these convergence problems (see Figure \ref{fig:007} ) consists in fixing the values ​​of the increments of plastic deformation and of its gradient to those of the previous time step.
%
These values thus become known contrary to the first version where they are unknowns of the problem.
%% 
%These values ​​thus become known contrary to the first version where they are unknowns of the problem.
%
%The resolution of the problem becomes similar to that of a purely elastic problem with initial deformations.
%
If the equilibrium equations are solved on the configuration at $ t$ and not at $ t + \Delta t$ there remain only a few weak nonlinearities linked to the objective derivatives of the stressess and moments.
At convergence we recover the converged values of the increments of plastic deformation and its gradient, which we distribute with the algorithm of projection on the criterion in elastic and plastic parts. 
%%
%%At convergence we recover the converged values ​​of the increments of plastic deformation and its gradient, which we distribute with the algorithm of projection on the criterion in elastic and plastic parts.
%
The projection algorithm is itself unchanged, the existence and uniqueness of the solution therefore remaining assured by the mean of the generalized standard character of the GLPD model at fixed porosity. 
%%
%%The projection algorithm is itself unchanged, the existence and uniqueness of the solution therefore remaining assured by the mean of the generalized standard character of the GLPD model at fixed porosity.
%
THe proposed modification is only valid for very small time steps.
%%
%%The proposed modification is only valid for very small time steps.
%
The error made on the increments of plastic deformation and its gradient is in $ \mathcal {O} ( (\Delta t)^2 ) $.
\newpage
%Indeed,
%%
%\begin{equation}
%\label{eqn:eq007}
%\Delta \boldsymbol{\varepsilon}^{i} = \Delta \boldsymbol{\varepsilon}^{e_i} + \Delta \boldsymbol{\varepsilon}^{p_i} \quad; \quad \Delta \boldsymbol{\varepsilon}^{i} = \Delta \boldsymbol{\varepsilon}^{e_i} + \Delta \boldsymbol{\varepsilon}^{p_{i-1}}; \nonumber
%\end{equation}
%%
%or $ \Delta \boldsymbol{\varepsilon}^{p_i} - \Delta \boldsymbol{\varepsilon}^{p_{i-1}}$ is proportional to $ (\Delta t)^2$.
%%
%For very small $ (\Delta t)^2$, $ \boldsymbol{\varepsilon}^{p_i} - \boldsymbol{\varepsilon}^{p_{i-1}} \rightarrow 0$.
%
%The agreement between the two curves is acceptable considering the experimental errors.
%
%The simulations on the TA30 specimen could not be completed.
%%
\section{Conclusion}
\label{sec:Con}
The contributions of this work can be summarized as follows:
\begin{itemize}
\item  We provide the exact consistent stiffness matrix for a porous materal model, the GLPD model in the framework of small deformations.
The expressions for the derivatives of the Cauchy stress tensors and the generalized moment stress tensors the model involved are derived. 
\item We have assessed the robustness of the formulation of the stiffness matrix proposed by comparing its numefical predictions with available experimental results of typical ductile fracture tests. 
The results show that quadratic convergence of the elasto-plastic ietrations was obtained by a slightly modify the original algorithm for the GLPD model.  
\end{itemize}


\newpage
\begin{thebibliography}{}
\bigskip
%%\bibitem[Collin et al.(2009)]{CCC09} Collin F., Caillerie D. and Chambon R. (2009). Analytical solutions for the thick-walled problem model with an isotropic elastic second gradient constitutive equation, {\it Int. J. of Solids Structures,} {\bf 46}, 3927-3937
%%\bibitem[Cosserat and Cosserat(1906)]{CC06} Cosserat E. and Cosserat F. (1906). Th\'eorie des corps d\'eformables, Hermann et fils Ed., Paris (in French).
%%\bibitem[Ahad et al.(2013)]{AESB93} Ahad, F., Enakoutsa, K., Solanki, K., and Bammann, D. (2013). Nonlocal Modeling in High Velocity Impact failure of Aluminum 6160-T6, {\it Int. J. Plast.,} in Press, doi:10.1016/j.ijplas.2013.10.001
%%\bibitem[Bergheau et al.(2013)]{BLP13} Bergheau, J-M, Leblond, J-B, Perrin, G. 2013. A New Numerical Implementation of a second-gradient model for plastic porous solids, with an application of ductile rupture tests, {\it Computer Methods in Applied Mechanics and Engineering,} {\bf 268}, 105-125    
\bibitem[Aravas(1987)]{A87} Aravas N. (1987). On the numerical integration of a class of pressure-dependent plasticity models,
  {\it Int. J. Num. Meth. Engng.}, {\bf 24}, 1395-1416.
%%  
\bibitem[Enakoutsa et al.(2007)]{ELP07} Enakoutsa K., Leblond J.B. and Perrin G. (2007). Numerical implementation and assessment of a phenomenological nonlocal model of ductile rupture, {\it Comput. Meth. Appl. Mech. Engng.}, {\bf 196}, 1946-1957
\bibitem[Enakoutsa(2007)]{E07}Enakoutsa K. (2007). Mod\'ele Non-locaux en rupture ductile des m\'etaux. Ph.D thesis, Universit\'e Pierre et Marie Curie (Paris VI) (in French).
\bibitem[Enakoutsa and Leblond(2009)]{EL09} Enakoutsa K., and Leblond J.B. (2009). Numerical implementation and assessment of the GLPD micromorphic model of ductile rupture, {\it Eur. J. Mech. A/Solids}, {\bf 28}, 445-460
%%\bibitem[Enakoutsa(2012)]{E12} Enakoutsa K. (2012). Some new applications of the GLPD micromorphic model for ductile fracture, {\it Math. Mech. Solids,} {\bf 19}(3), 242-259
%\bibitem[Enakoutsa et al.(2012a)]{ESAB12} Enakoutsa, K., Solanki, K., Ahad, F., Tjiptowidjojo, Y., and Bammann, D. (2012). Using Damage Delocalization to Model Localization Phenomena in Bammann-Chiesa-Johnson Metals, {\it J. Eng. Mater. Tech.,} 134(4)
%\bibitem[Enakoutsa et al.(2012b)]{ESATB12} Enakoutsa, K., Solanki, K., Ahad, F., Tjiptowidjojo, Y., and Bammann, D. (2012). Damage smoothing effects in a delocalized rate sensitivity model for metals, {\it Theoretical and Applied Mechanical Letters,} {\bf 2}(5): 5-051005 
\bibitem[Enakoutsa(2012)]{E12} Enakoutsa, K., 2012.``Some new Applications of the GLPD Micromorphic Model of Ductile Fracture,'' {\it Mathematics and Mechanics of Solids,} {\bf 19}(3), 242-259
%%
\bibitem[Enakoutsa (2013a)]{E13}Enakoutsa, K. (2013). Exact results for the problem of a hollow sphere subjected to hydrostatic tension and made of micromorphic plastic porous materials, {\it Mech. Res. Commun.,} {\bf 49}, 1-7
%\bibitem[Enakoutsa(2013b)]{E13b} Enakoutsa, K. (2013). An Analytic Benchmark Solution to the Problem of a Generalized Plane Strain Hollow Cylinder made of Micromorphic Plastic Porous Metal and Subjected Axisymmetric Loading Conditions, {\it Math. Mech. Solids,} In press 
%%\bibitem[Enakoutsa (2013)]{E13b} Enakoutsa K. (2013). A Constitutive Generalized Elasticity Law for a Micromorphic Second Gradient Elastic Material, {\it Theoretical and Applied Mechanics Letters,} submitted 
%%\bibitem[Gao(2003a)]{G03a} Gao X.L. (2004). Elasto-plastic analysis of an internally pressurized thick-walled cylinder using a strain gradient plasticity theory, {\it Int. J. Solids Structures,} {\bf 40}, 6445-6455
%%\bibitem[Gao and Park(2006)]{GP06} Gao X.L. and Park S.K. (2006). Supplemental notes on derivations (81 pages, available upon request)  
%%\bibitem[Gao and Park(2007)]{GP07} Gao X.L. and Park S.K. (2007). Variational formulation of a simplified strain gradient elasticity theory and its application to a pressurized thick-walled cylinder problem, {\it Int. J. Solids Structures,} {\bf 44}, 7486-7499 
%%\bibitem[Gao(2003b)]{G03b} Gao X.L. (2003). Strain gradient plasticity solution for an internally pressurized thick-walled spherical shell of an elasto-plastic material, {\it Mech. Res. Commun.,} {\bf 30}, 411�420  
%%\bibitem[Germain(1973b)]{Ger73b} Germain P. (1973). The Method of virtual power in continuum mechanics. Part 2: Microstructure, {\it SIAM J. of Appl. Math.}, {\bf 25}, 556-575
%%\bibitem[Van der Giessen and Needleman(1995)]{GN95} Van der Giessen E. and Needleman A. (1995). Discrete dislocation plasticity: a simple planar model, {\it Modell. and Simul. Mater. Sci. Engng.,} {\bf 3}, 689-
\bibitem[Forest(1998)]{F98} S. Forest, 1998. ``Mechanics of generalized continua: construction by homogenization,'' \textit{J. Phys. IV,} {\bf 8}, 39-48.
%%
\bibitem[Forest {\it et al}.(2000)]{FBC00} Forest S., Barbe F. and Cailletaud G. (2000). Cosserat modelling of size effects in the mechanical behaviour
of polycrystals and multi-phase materials, {\it Int. J. Solids Structures}, {\bf 37}, 7105-7126.
\bibitem[Gologanu et al.(1997)]{GLPD97} Gologanu M., Leblond J.B., Perrin G. and Devaux J. (1997). Recent extensions of Gurson's model for porous ductile metals, in: {\it Continuum Micromechanics}, CISM Courses and Lectures 377, P. Suquet ed., Springer, pp. 61-130
\bibitem[Gurson(1977)]{G77} Gurson A.L. (1977). Continuum theory of ductile rupture by void nucleation and growth: Part I - yield criteria and flow rules for porous ductile media, {\it ASME J. Engng. Mater. Technol.}, {\bf 99}, 2-15
%%\bibitem[Hall and Petch(1950)]{HP50} Hall E.O. and Petch N.J., 1950. The cleavage strength of crystals, {\it J. Iron Steel Institute,} {\bf 174}, 25-28
\bibitem[Halphen and Nguyen(1975)]{HN75} Halphen B. and Nguyen Q.S. (1975). Sur les mat\'eriaux standards g\'en\'eralis\'es, {\it Journal de M\'ecanique},
%%
\bibitem[Hill(1967)]{H67} Hill, R., 1967. The essential structure of constitutive laws for metal composites and polycrystals, {\it Journal of Mechanics and Physics of Solids}, {\bf 15}, 79-95
%%
\bibitem[Leblond et al.(1994)]{Leb94} Leblond, J.B., Perrin, G., and Devaux, J.,(1994). Bifurcation Effects in Ductile Metals with Nonlocal Damage, {\it ASME J. Applied . Mech.}, {\textbf 61}, 236-242.
%\bibitem[Nahshon and Hutchinson(2008)]{NH08} Nahshon, K., Hutchinson, J.W., 2008. Modification of the Gurson model for shear failure, {\it Eur. J. Mech. A/Solids,} {\bf 27}, 1-17.
\bibitem[Mandel(1964)]{M64} Mandel, J., 1964. Contribution th\'eorique \`a l'\'etude de l'\'ecrouissage et des lois d'\'ecoulement plastique, {\it Proceedings of the 11th International Congress on Applied Mechanics}, Springer, pp. 502-509 (in French)
\bibitem[Matsushima {\it et al}.(2000)]{MCC00} Matsushima T., Chambon R. and Caillerie D. (2000). Second gradient models as a particular case of microstructured models:
  a large strain finite element analysis, {\it Comptes-Rendus Acad. Sc. Paris S\'erie IIb}, {\bf 328}, 179-186
%%
\bibitem[Nguyen(1977)]{N77} Nguyen Q.S. (1977). On the elastic plastic initial-boundary value problem and its numerical integration, {\it Int. J. Numer.
  Meth. Engng.}, {\bf 11}, 817-832. 
\bibitem[Perrin and Leblond(1990)]{PL90} Perrin, G, and Leblond, J-B, 1990. Analytical study of a hollow sphere made of plastic porous material and subjected to hydrostatic tension: Application to some problems in ductile fracture of metals, 
\bibitem[Perrin and Leblond(2000)]{PL00} Perrin, G, and Leblond, JB (2000). Accelerated void growth in porous ductile solids containing two populations of cavities. {\it Int. J. Plast.,} {\bf 16}, 91-120
{\it Int. J. Plast.,} {\bf 6}:677-699
\bibitem[Pijaudier-Cabot and Bazant(1987)]{PijB87} Pijaudier-Cabot, G. and Bazant, Z.P., (1987). Nonlocal Damage Theory, {\it ASCE J. Engrg. Mech.}, {\textbf 113}, 1512-1533
%%\bibitem[Mindlin(1964)]{Min64} Mindlin R.D. (1964). Microstructure in linear elasticity, {\it Arch. Rational Mech. Anal.}, {\bf 12}, 51-78
%%\bibitem[Mindlin(1965)]{Min65} Mindlin R.D. (1965). Second gradient of strain and surface-tension in linear elasticity, {\it Int. J. Solids Structures}, {\bf 1}, 417-738
%\bibitem[Rice and Tracey(1969)]{RT69} Rice, J. R. and Tracey, D. M. (1969). On the Ductile Enlargement of Voids in Triaxial Growth Holes. {\it J. Mech. Phys. Solids,} {\bf 17}, 210-217  
%\bibitem[Tvergaard(1981)]{T81} Tvergaard V. (1981). Influence of voids on shear band instabilities under plane strain conditions, {\it Int. J. Fracture}, {\bf 17}, 389-407
%
\bibitem[Rousselier and Murdy (1981)]{RM81} Rousselier, G and Mudry F, 1981. Etude de la Rupture Ductile de l'Acier Faiblement Allie en Mn-Ni-Mo pour
Cuves de Reacteurs a Eau Ordinaire sous Pression, Approvisionne sous la forme d'une Debouchure de Tubulure. Resultats du
Programme Experimental, EdF Centre des Renardieres Internal Report HT/PV D529 MAT/T43
\bibitem[Shu {\it et al}.(1999)]{SKF99} Shu J., King W. and Fleck N. (1999). Finite elements for materials with strain gradient effects, {\it Int. J. Numer.
  Methods Engng.}, {\bf 44}, 373-391.
\bibitem[Tvergaard(1981)]{T81}\textcolor{blue} {Tvergaard V. (1981). Influence of voids on shear band instabilities under plane strain conditions, {\it Int. J. Fracture}, {\bf 17}, 389-407.}
\bibitem[Tvergaard and Needleman(1984)]{TN84} Tvergaard V. and Needleman A., 1984. Analysis of cup-cone fracture in a round tensile bar, {\it Acta Metall.}, {\bf 32}, 157-169
\bibitem[Tvergaard and Needleman(1997)]{TN97} Tvergaard V. and Needleman A.(1997). Nonlocal effects on localization in a void-sheet, {\it Int. J. Solids Structures,} {\bf 34}, 2221-2238
\bibitem[Tvergaard and Needleman(1995)]{TN95} Tvergaard V. and Needleman A. (1995). Effects of nonlocal damage in porous plastic solids, {Int. J. Solids Structures,} {\bf 32}, 1063-1077 s
%\bibitem[Aifantis(1984)]{Aif84} Aifantis, E.C., 1984. ``On the Microstructural Origin of Certain Inelastic Models,'' \textit{Journal of Engng Mater. Tech.} {\bf 106}, 326-334
%%
%%
%\bibitem[Tvergaard and Needleman(1984)]{TN84} \textcolor{blue}{Tvergaard V. and Needleman A. (1984). Analysis of cup-cone fracture in a round tensile bar, {\it Acta Metall.},{\bf 32}, 157-169.}
%%
%%\bibitem[Tvergaard and Needleman(1984)]{TN84} \textcolor{blue}{Tvergaard V. and Needleman A. (1984). Analysis of cup-cone fracture in a round tensile bar, {\it Acta Metall.},{\bf 32}, 157-169.}
%%
%\bibitem[Tvergaard(1989)]{TVAIAM83}{\textcolor{blue}{Tvergaard, V. (1989). Material failure by void growth to coalescence. In Advances in applied Mechanics (Vol. 27, pp. 83-151).Elsevier.}}
%%
%\bibitem[Tvergaard(1996)]{TV83}{\textcolor{blue}{Tvergaard, V. and Needleman, A. (1996). Nonlocal effects of localization in a void-sheet. International Journal of solid structures (Vol. 34, No. 18, pp, 2221-2338). Elsevier.}}.
%%\bibitem[Aifantis(1995)]{Aif95} Aifantis, E.C., 1995.``From Micro-Plasticity To Macro-Plasticity: The Scale-Invariance Approach,'' {\it J. of Eng. Materials and Technology-Transaction of the ASME}, {\bf 117}(4), 352-355
%%
%\bibitem[Bammann and Aifantis(1982)]{BA82} Bammann, D.J. and Aifantis, E.C., 1982. On a Proposal for Continuum with Microstructure. Acta Mechanica, 45(1-2), 91-121
%%
%%\bibitem[Bammann(2000)]{B00}Bammann, J.D., 2000. ``A model of crystal plasticity containing a natural length scale,'' {\it Mat. Sci. Eng.,} A309-310, 406-410
%%\bibitem[Bazant and Pijaudier-Cabot(1989)]{BPC89} Bazant, Z. and Pijaudier-Cabot, G., 1989. ``Measurement of Characteristic Length of Nonlocal Continuum,'' {\it J. Eng. Mech.,} {\bf 115}(4), 755-767
%%\bibitem[Bell(1979-81)]{B79} Bell, J.F., 1979-81. ``A Physical basis for continuum theories of finite plasticity, I$&$II,'' {\it Archive for Rational Mechanics and Analysis,} Part I {\bf 70}, 319-(1979); Part II {\bf 75}, 103-(1981) 
%%\bibitem[Bell and Khan(1980)]{BK80} Bell, J.F., Khan, A.S., 1980. ``Finite plastic strain associated annealed copper during non-proportional loading,'' {\it Int. J. Solids Structures,} {\bf 16}, 683-
%%\bibitem[Bergheau et al.(2013)]{BLP13} Bergheau, J-M, Leblond, J-B, and Perrin, G. 2013. ``A New numerical implementation of a second-gradient model for plastic porous solids, with an application of ductile rupture tests,'' {\it Computer Methods in Applied Mechanics and Engineering,} {\bf 268}, 105-125
%%
%\bibitem[Dell'Isola et al.(2009)]{DSV09} F. Dell'Isola, G. Sciarra G., S. Vidoli, 2009. ``Generalized Hooke's Law for Isotropic Second Gradient materials,'' {\it Proc. R. So. A} {\textbf 465}(\textbf{no. 2107}), 2177-2196
%%
%%\bibitem[Clayton et al.(2006)]{CMB06} Clayton, J.D., McDowell, D.L., and Bammann, D.J., 2006. ``Modeling dislocations and disclinations with finite micropolar elastoplasticity,'' {\it International Journal of Plasticity,} {\bf 22,} 210-256.
%%
%\bibitem[Collin et al.(2009)]{CCC09} Collin, F., Caillerie, D., Chambon, R., 2009. ``Analytical solutions for the thick-walled cylinder problem model with an isotropic elastic second gradient constitutive equation,'' {\it Int. J. Solids and Structures,} {\bf 46}, 3927-3937
%%
%\bibitem[Cosserat and Cosserat(1906)]{CC06} Cosserat, E. and Cosserat, F., 1906. ``Th\'eorie des corps  d\'eformables,'' Hermann et fils Ed., Paris (in French)
%%
%\bibitem[Enakoutsa(2007)]{E07} Enakoutsa, K., 2007. Mod\`ele Non-locaux en Rupture Ductile des M\'etaux. Ph.D thesis, Universit\'e Pierre et Marie Curie (Paris VI) (in French)
%%
%\bibitem[Enakoutsa and Leblond(2009)]{EL09} Enakoutsa K., and Leblond J.B., 2009.``Numerical Implementation and Assessment of the GLPD Micromorphic Model of Ductile Rupture,'' {\it European Journal of Mechanics A/Solids}, {\bf 28}, 445-460 
%%
%%
%\bibitem[Enakoutsa(2013)]{E13}Enakoutsa, K. 2013 ``The method of virtual power in a micromorphic theory of ductile fracture in metals,'' International Journal of Applied Mathematics and Multiscale Mechanics, 2(4), 311322 

%\bibitem[Enakoutsa(2015)]{E15} Enakoutsa, K., 2015 ``Analytical applications and effective properties of a second gradient isotropic elastic material model,'' Zeitschrift fr angewandte Mathematik und Physik, 66(3), 1277-1293
%%
%%\bibitem[Enakoutsa (2013a)]{E13}Enakoutsa, K., 2013. ``Exact results for the problem of a hollow sphere subjected to hydrostatic tension and Made of micromorphic plastic porous materials,'' {\it Mechanics Research Communications,} {\bf 49}, 1-7
%%
%\bibitem[Enakoutsa(2013)]{E13} Enakoutsa, K., 2013. ``An Analytic Benchmark Solution to the Problem of a Generalized Plane Strain Hollow Cylinder made of Micromorphic Plastic Porous Metal and Subjected Axisymmetric Loading Conditions,'' {\it Journal of Mathematics and Mechanics of Solids,} DOI: 10.1177/1081286513513457
%%
%\bibitem[Efremidis et al.(2004)]{ERA04} Efremidis, G., Rambert, G., Aifantis, E.C., 2004. ``Gradient elasticity and size effects for pressurized thick hollow cylinder,'' \textit{Journal of the Mechanical Behavior of Materials,} {\bf 15}(3), 169-184
%%
%\bibitem[Enakoutsa(2012b)]{E2012} Enakoutsa, K., 2012.``Modeling ductile fracture in metals involving two populations of voids�Cinfluence of continuous nucleation of secondary voids upon growth and coalescence of primary voids,'' {\it Mathematics and Mechanics of Solids,} {\bf 18}(3), 323-345
%\bibitem[Enakoutsa et al.(2007)]{ELP07} Enakoutsa K., Leblond J.B. and Perrin G., 2007.``Numerical Implementation and Assessment of a Phenomenological Nonlocal Model of Ductile Rupture,'' {\it Computational Methods in Applied Mechanics and Engineering}, {\bf 196}, 1946-1957
%\bibitem[Eringen and Suhubi(1964)]{ES64} Eringen, A. C. and Suhubi, E. S., 1964. ``Nonlinear theory of simple micro-elastic solids,'' {\it I. Int. J. Engng Sci.} {\bf 2,} 189-203.
%\bibitem[Eringen(1992)]{E92} Eringen, A. C., 1992. Vistas of nonlocal continuum physics. Int. J. Engng Sci. 30(10), 1551-1565
%\bibitem[Gao(2003a)]{G03a} Gao X.L. (2004). Elasto-plastic analysis of an internally pressurized thick-walled cylinder using a strain gradient plasticity theory, {\it Int. J. Solids Structures,} {\bf 40}, 6445-6455
%%\bibitem[Gao and Park(2006)]{GP06} Gao X.L. and Park S.K.(2006). Supplemental notes on derivations (81 pages, available upon request)
%%
%%\bibitem[Gao(2003b)]{G03b} Gao X.L.(2003). Strain gradient plasticity solution for an internally pressurized thick-walled spherical shell of an elasto-plastic material, {\it Mech. Res. Commun.,} {\bf 30}, 411??420
%%
%%
%%
%\bibitem[Gao et al.(2008)]{GPM08} Gao X.L., Park S.K., Ma, H. M., 2008. ``Analytical Solution for a Pressurized Thick-Walled Spherical Shell Based on a Simplied Strain Gradient Elasticity Theory,'' {\it Mathematics and Mechanics of Solids,} {\bf 14}, 747??758
%%
%\bibitem[Gao et al.(2007)]{GP07} X.L Gao, and S. Park, 2007. ``Variational formulation of a simplified strain gradient elasticity theory and its application to a pressurized thick-walled cylinder problem.'' \textit{Int. J. Solids Struct.,} \textbf{44,} 7486-7499.    
%%
%\bibitem[Germain(1973)]{Ger73} Germain, P., 1973. ``The Method of Virtual Power in Continuum Mechanics. Part 2: Microstructure,'' {\it SIAM Journal of Applied Mathematics}, \textbf {25}, 556-575
%%
%\bibitem[Germain(1986)]{G86} Germain, P., \textbf{M\'ecanique} (Tome II), Ecole Polytechnique (1986)
%% 
%%\bibitem[Van der Giessen and Needleman(1995)]{GN95} Van der Giessen, E. and Needleman, A., 1995. ``Discrete dislocation plasticity: a simple planar model,''{\it Modelling Simul. Mater. Sci. Eng.}, \textbf{3}, 689-
%\bibitem[Gao(2003a)]{G03a} Gao, X.L., 2004.``Elasto-plastic analysis of an Internally Pressurized Thick-walled Cylinder using a Strain Gradient Plasticity Theory,'' {\it Int. J. Solids and Structures,} {\bf 40}, 6445-6455
%%
%\bibitem[Gologanu et al.(1997)]{GLPD97} Gologanu M., Leblond J.B., Perrin G. and Devaux J., 1997. ``Recent extensions of Gurson's model for porous ductile metals,'' in: {\it Continuum Micromechanics}, CISM Courses and Lectures 377, P. Suquet ed., Springer, pp. 61-130
%%
%\bibitem[Hall and Petch(1950)]{HP50} Hall, E.O. and Petch, N.J., 1950.``The Cleavage Strength of Crystals,'' {\it J. Iron and Steel Inst.,} \textbf{174}, 25-28
%\bibitem[Fleck and Hutchinson (1993)]{FH93} Fleck, N. A., and Hutchinson, J.W., 1993. ``A Phenomenological Theory for Strain Gradient Effects in Plasticity,'' {\it Journal of the Mechanics and Physics of Solids,} {\bf 41,} 1825-1857
%%
%\bibitem[Fleck and Hutchinson (1997)]{FH97} Fleck, N. A. and Hutchinson, J.W., 1997. ``Strain Gradient Plasticity,'' {\it Adv. Appl. Mech.,} {\bf 33,} 295-361
%\bibitem[Fleck et al.(1994)]{FMAH94} Fleck N.A., Muller G.M., Ashby M.F., and Hutchinson J.W., 1994. ``Strain Gradient Plasticity: Theory and Experiment,'' {\it Acts Metallurgic et Materialia,} {\bf 42}(2), 475-487
%\bibitem[Gurson(1977)]{G77} Gurson A.L. 1977.``Continuum Theory of Ductile Rupture by Void Nucleation and Growth: Part I - Yield Criteria and Flow Rules for Porous Ductile Media,'' {\it ASME J. Engng. Materials Technol.}, \textbf{99}, 2-15
%\bibitem[McMeeking(1982)]{M82} McMeeking, R.M., 1982. ``The Finite Strain Tension Torsion Test of a Thin-Walled Tube of Elastic-Plastic Material,'' {\it Int. J. Solids Struct.,} {\bf 18}, 199-204
%%
%\bibitem[Kiang et al.(1998)]{KEADD98} C.-H. Kiang, M. Endo, P. M. Ajayan, G. Dresselhaus, and M. S. Dresselhaus, 1998. ``Size Effects in Carbon Nanotubes,'' \textit{Physics Review Letters} {\bf 81}, 1869- 
%%
%\bibitem[Koiter(1964)]{K64} Koiter, W.T., 1964. Couple stress in the theory of elasticity I-II. Proc. Nederl. Akad. Wetensch. 67, 17-44
%%
%\bibitem[Mindlin(1964)]{Min64} Mindlin, R.D., 1964. ``Microstructure in Linear Elasticity,'' {\it Archive of Rational Mechanics and Analysis}, \textbf{12}, 51-78
%%
%\bibitem[Mindlin(1965)]{Min65} Mindlin R.D., 1965. ``Second Gradient of Strain and Surface-tension in Linear Elasticity,'' {\it International Journal of Solids and Structures}, \textbf{1}, 417-738%
%\bibitem[McLachlan(1941)]{M41} McLachlan, N.W., 1941. Bessel Functions for Engineers. Oxford University Press, London.
%%
%\bibitem[Maugin(2011)]{Maug11} Maugin, G. A., 2011. A historical perspective of generalized continuum mechanics. In H. Altenbach, G.A. Maugin and V. Erofeev, editors, Mechanics of generalized continua - From micromechanical basis to engineering applications, 1-17. Springer, Berlin
%%
%\bibitem[Mindlin and Tiersten(1962)]{MT62} Mindlin, R.D. and Tiersten, H.F., 1962. Effects of couple stresses in linear elasticity. Arch. Rat. Mech. Anal. 11, 415-448
%%
%\bibitem[Muhlich et al.(2012)]{MZK12} Muhlich, U., Zybell, L., Kuna, M., 2012. ``Estimation of material properties for linear elastic gradient effective media,'' \textit{Eur. J. Mech. A-Solid}, {\bf 31}(1) (2012), 117-
%%
%\bibitem[Abu Al-Rub and Faruk(2012)]{AAF12} Abu Al-Rub, R.K. and Faruk, A.N.M., 2012. ``Prediction of micro and nano indentation size effects from spherical indenters,'' {\it Mechanics of Advanced Materials and Structures,} \textbf{19}1-3, 119-128
%\bibitem [Pineau and Joly(1991)]{PJ91} Pineau, A. and Joly, P., 1991. ``Local versus global approaches of elastic-plastic fracture mechanics. Application to ferritic steels and a cast duplex stainless steel,'' in {\em Effect Assessment in Components - Fundamentals and Applications}, J.G. Blauel and K.H. Schwalbe, eds., ESIS, European Group on Fracture Publication 9.  
%\bibitem[Solanki and Bammann(2010)]{SB10} Solanki, K. and Bammann, D.J., 2010. A Thermodynamic% Framework for a Gradient Theory of Continuum Damage. Acta Mechanica, 213, 271-738
%\bibitem[Taylor and Quinney(1931)]{TQ31} Taylor, G.I., Quinney, H., 1931. The Plastic Distorsion of Metals, {\it Phil. Trans. Roy. Soc. Lon.}, {\bf 230}(A), 323-362
%\bibitem[Toupin(1962)]{T62} Toupin, R. A., 1962. ``Elastic materials with couple-stresses,'' {\it Arch. Rational Mech. Anal.,} {\bf 11,} 385-414 
%\bibitem[Tvergaard(1981)]{T81} Tvergaard, V., 1981.``Influence of Voids on Shear Band Instabilities under Plane Strain Conditions,'' {\em International Journal of Fracture}, {\bf 17}, 389-407
%\bibitem[Tvergaard and Needleman(1984)]{TN84} Tvergaard, V. and Needleman, A., 1984.``Analysis of Cup-cone Fracture in a Round Tensile Bar,'' {\em Acta Metallurgica}, {\bf 32}, 157-169
%\bibitem[Tvergaard and Needleman(1997)]{TN97} Tvergaard, V. and Needleman, A., 1997.``Nonlocal Effects on Localization in a Void-sheet'', {\it International Journal of Solids and Structures,} {\bf 34}, 2221-2238
%\bibitem[Tvergaard and Needleman(1995)]{TN95} Tvergaard V. and Needleman A., 1995.``Effects of Nonlocal Damage in Porous Plastic Solids'', {\it Int. J. Solids Struct.,} {\bf 32}, 1063-1077
%\bibitem[Chu and Needleman(1980)]{CN80} Chu, C.C. and Needleman, A., 1980. ``Void nucleation effects on shear localization in porous plastic solids,'' {\em ASME J. Engng. Materials Technol.}, {\bf 102}, pp. 249-256.     
%\bibitem[Tvergaard(1981)]{T81} Tvergaard V., 1981. ``Influence of voids on shear band instabilities under plane strain conditions,'' {\it Int. J. Fracture}, {\bf 17}, 389-407.
%\bibitem[Tvergaard and Needleman(1984)]{TN84} Tvergaard V. and Needleman A. (1984). Analysis of cup-cone fracture in a round tensile bar, {\it Acta Metall.}, {\bf 32}, 157-169
%\bibitem[Pijaudier-Cabot and Bazant(1987)]{PCB87} Pijaudier-Cabot G. and Bazant Z.P. (1987). Nonlocal damage theory, {\it ASCE J. Engng. Mech.}, {\bf 113}, 1512-1533.
%\bibitem[Leblond et al.(1994)]{LPD94} Leblond J.B., Perrin G. and Devaux J., 1994. ``Bifurcation effects in ductile metals with nonlocal damage,'' {\it ASME J. Appl. Mech.}, {\bf 61}, 236-242.
%\bibitem[Rousselier(1987)]{R87} Rousselier, G., 1987. ``Ductile fracture models and their potential in local approach of fracture,'' {\it Nuclear Engrg. Design,} {\bf 105}, 97�C111
%\bibtem[Howard et al.(1994)]{HLB94} Howard, I.C., Li, Z.H., Bilby, B.A., 1994. ``Ductile crack growth predictions for large center cracked panels by damage modeling using 3-D ??nite element analysis,'' {\it Fatigue and Fracture Engrg. Mater. Struct.} {\bf 17}, 959�C969
%\bibitem[Xia et al.(1995)]{XSH95} Xia, L., Shih, C.F., Hutchinson, J.W., 1995. ``A computational approach to ductile crack growth under large scale yielding conditions,'' {\it J. Mech. Phys. Solids,} {\bf 43}, 398�C413.
%\bibitem[Tvergaard(1990)]{T90} Tvergaard, V., 1990. ``Material failure by void growth,'' {\it Adv. Appl. Mech.} {\bf 27}, 83�C151
%\bibitem[Nahshon and Hutchinson(2008)]{NH08} K. Nahshon, J.W. Hutchinson, 2008. ``Modi??cation of the Gurson Model for shear failure,'' {\it European Journal of Mechanics A/Solids,} {\bf 27,} 1�C17
%\bibitem[Barsoum and Faleskog]{BF07} Barsoum, I., Faleskog, J., 2007. ``Rupture in combined tension and shear: Experiments,'' {\it Int. J. Solids Structures,} {\bf 44}, 1768�C1786
%\bibitem[Sokolowski(1970)]{S70} Sokolowski, M., 1970. Theory of Couple-Stresses in Bodies with Constrained Rotations. CISM courses and lectures 26, Berlin, Germany: Springer
%%
%%
%\bibitem[Zhao(2011)]{Z11} Zhao, J., 2011. ``A unifed theory for cavity expansion in cohesive-frictional micromorphic media,'' \textit{International Journal of Solids and Structures,} {\bf 48}, 1370-1381
%%  
%%
%%\bibitem{Bonora(1997)}{\textcolor{blue}{Bonora, N. (1997). A nonlinear CDM model for ductile failure. Engineering Fracture Mechanics, 58(1), 11-28}}
%%
%\bibitem[Enakoutsa {\it et al}.(2007)]{ELP07} \textcolor{blue}{Enakoutsa K., Leblond J.B. and Perrin G. (2007). Numerical implementation and assessment of a phenomenological nonlocal model of ductile rupture, {\it Comput. Meth. Appl. Mech. Engng.}, {\bf 196}, 1946-1957.}
%%
%\bibitem[Enakoutsa(2013)]{Koffi2013}\textcolor{blue}{Enakoutsa. "Exact Results for the Problem of a hollow sphere subjected to hydrostatic tension and made of micromorphic porous materials." Mechanics Research Communications:(2013)}
%%
%\bibitem[Krishnan(2018)]{JIHES92}{\textcolor{blue}{Krishnan, S., Sasikala, G., Moitra, A., Albert, S., \& Bhaduri, A. (2018). Ductile damage parameters and far field J-integral for high hardening steel. International Journal of Structural Integrity, 9(2), 153-167.}}
%%
%\bibitem [Komori(2017)]{PE207}{\textcolor{blue}{Komori, K. (2017). Improvement and validation of the ellipsoidal void model for predicting ductile fracture. Procedia engineering, 207, 2036-2041}}
%%
%\bibitem[Gologanu {\it et al}.(1997)]{GLPD97}{\textcolor{blue}{Gologanu M., Leblond J.B., Perrin G. and Devaux J. (1997). Recent extensions of Gurson's model for porous ductile metals, in: {\it Continuum Micromechanics}, P. Suquet, ed., CISM Courses and Lectures No. 377, Springer, pp. 61-130.}}
%%
%%
%\bibitem[Pijaudier-Cabot and Bazant(1987)]{PCB87} \textcolor{blue}{Pijaudier-Cabot G. and Bazant Z.P. (1987). Nonlocal damage theory, {\it ASCE J. Engng. Mech.}, {\bf 113}, 1512-1533.}
%%
%\bibitem[Sadasue(2004)]{ASME1669}{\textcolor{blue}{Sadasue, T., Igi, S., Kubo, T., Ishikawa, N., Endo, S., Glover, A., ... \& Toyoda, M. (2004, January). Ductile cracking evaluation of $X80/X100$ high strength linepipes. In 2004 International Pipeline Conference (pp. 1661-1669). American Society of Mechanical Engineers.}}
%%
\end{thebibliography}
\end{document}